\documentclass[12pt,twoside]{amsart}
\usepackage{amssymb}

\usepackage[all]{xy}
\nonstopmode

\numberwithin{equation}{section} 

\font\tengothic=eufm10 scaled\magstep 1 \font\sevengothic=eufm7
scaled\magstep 1
\newfam\gothicfam
       \textfont\gothicfam=\tengothic
       \scriptfont\gothicfam=\sevengothic
\def\goth#1{{\fam\gothicfam #1}}

\newtheorem{theorem}{Theorem}[section]
\newtheorem{lemma}[theorem]{Lemma}

\newtheorem{corollary}[theorem]{Corollary}
\newtheorem{conjecture}[theorem]{Conjecture}

\theoremstyle{definition}
\newtheorem{definition}[theorem]{Definition} 
\newtheorem{remark}[theorem]{Remark}
\newtheorem{example}[theorem]{Example}

\newcommand{\pd}{\operatorname{pd}}
\newcommand{\codim}{\operatorname{codim}}
\newcommand{\coker}{\operatorname{coker}}

\newcommand{\Supp}{\operatorname{Supp}}

\newcommand{\Hom}{\operatorname{Hom}}

\newcommand{\ext}{\operatorname{ext}}
\newcommand{\Ext}{\operatorname{Ext}}

\newcommand{\depth}{\operatorname{depth}}

\newcommand{\im}{\operatorname{im}}

\newcommand{\Hi}{\operatorname{Hilb}}
\DeclareMathOperator{\GradAlg}{GradAlg}

\newcommand{\Proj}{\operatorname{Proj}}
\newcommand{\Spec}{\operatorname{Spec}}
\newcommand{\proj}[1]
{ \mathchoice
            { {\mathbb P}^{#1} }
            { {\mathbb P}^{#1} }
            { {\mathbb P}^{#1} }
            { {\mathbb P}^{#1} }
          }

\newcommand{\Tor}{\operatorname{Tor}}

\newcommand{\cA}{{\mathcal A}}
\newcommand{\cB}{{\mathcal B}}
\newcommand{\cH}{{\mathcal H}}
\newcommand{\cI}{{\mathcal I}}

\newcommand{\cF}{{\mathcal F}}
\newcommand{\cG}{{\mathcal G}}
\newcommand{\cO}{{\mathcal O}}

\newcommand {\PP}{\mathbb{P}}
\newcommand {\AAA}{\mathbb{A}}

\newcommand {\VV}{\mathbb{V}}
\newcommand {\ra}{\longrightarrow}

\begin{document}
\title[]{Deformations of modules of maximal grade and the Hilbert scheme at
  determinantal schemes.}

\author[Jan O.\ Kleppe]{Jan O.\ Kleppe} 
 \address{ {\SMALL Faculty of Engineering,
         Oslo University College,
         Pb. 4 St. Olavs plass, N-0130 Oslo,
         Norway } }
\email{JanOddvar.Kleppe@iu.hio.no}


\date{\today}



\begin{abstract} 
  Let $R$ be a polynomial ring and $M$ a finitely generated graded $R$-module
  of maximal grade (which means that the ideal $ I_t(\cA)$ generated by the
  maximal minors of a homogeneous presentation matrix, $\cA$, of $M$ has
  maximal codimension in $R$). Suppose $X:=\Proj(R/I_t(\cA))$ is smooth in a
  sufficiently large open subset and $\dim X \ge 1$. Then we prove that the
  local graded deformation functor of $M$ is isomorphic to the local Hilbert
  (scheme) functor at $X \subset \Proj(R)$ under a week assumption 
  which holds if $\dim X \ge 2$. Under this assumptions we get that the
  Hilbert scheme is smooth at $(X)$, and we give an explicit formula for the
  dimension of its local ring. As a corollary we prove a conjecture of R.
  \!M.\! Mir\'o-Roig and the author that the closure of the locus of standard
  determinantal schemes with fixed degrees of the entries in a presentation
  matrix is a generically smooth component $V$ of the Hilbert scheme. Also
  their conjecture on the dimension of $V$ is proved for $\dim X \ge 1$. The
  cohomology $H^i_{*}({\mathcal N}_X)$ of the normal sheaf of $X$ in $
  \Proj(R)$ is shown to vanish for $1 \le i \le \dim X-2$. Finally the
  mentioned results, slightly adapted, remain true replacing $R$ by any
  Cohen-Macaulay quotient of a polynomial ring. 
\end{abstract}


\maketitle



\section{Introduction} \label{intro} Determinantal objects are central in many
areas of mathematics. In algebraic geometry determinantal 
schemes defined by the vanishing of the $p \times p$\,-minors of a homogeneous
polynomial matrix, may be used to describe classical schemes such as rational
normal scrolls and other fibered schemes, Veronese and Segre varieties and
Secant schemes to rational normal curves and Segre varieties (\cite{Ha},
\cite{BF}). Throughout the years many nice properties are detected for
determinantal schemes, 
e.g. they are arithmetically Cohen-Macaulay with rather well understood free
resolutions and singular loci, see \cite{e-h}, \cite{e-n}, \cite{dan2},
\cite{wey}, and see \cite{b-v}, \cite{BH}, \cite{eise}, \cite{Go},
\cite{KMMNP}, \cite{dan}, \cite{M} for history and other important
contributions. 

In this paper we study the Hilbert scheme along the locus of determinantal
schemes. More precisely we study deformations of modules of maximal grade over
a polynomial ring $R$ and establish a very strong connection to corresponding
deformations of determinantal schemes in $ \PP^{n}$. Recall that the grade $g$
of a finitely generated graded $R$-module $M$ is the grade of its annihilator
$I:= {\rm ann}(M)$, i.e. $g=\depth_{I} R =\dim R- \dim R/I$. We say a scheme
$X\subset \PP^{n}$ of codimension $c$ is {\it standard determinantal} if its
homogeneous saturated ideal is equal to the ideal $I_t(\cA)$ generated by the
$t \times t$ minors of some homogeneous $t \times (t+c-1)$ matrix
$\cA=(f_{ij})$, $f_{ij} \in R$. If $M$ is the cokernel of the map determined
by $\cA$, then $g=c$ because the radicals of $I$ and $I_t(\cA)$ are equal.
Moreover $M$ has maximal grade if and only if $X=\Proj(A), A:=R/I_t(\cA)$ is
standard determinantal. In this case $ {\rm ann}(M)=I_t(\cA)$ for $c \ge 2$ by
\cite{eis}. 
                                         

Let $ \Hi ^p(\PP^{n})$ be the Hilbert scheme parameterizing closed subschemes
of $\PP^{n}$ of dimension $n-c \ge 0$ and with Hilbert polynomial $p$. Given
integers $a_0\le a_1\le ...\le a_{t+c-2}$ and $b_1\le ...\le b_t$, $t \ge 2$,
$c \ge 2$, we denote by $W_s(\underline{b};\underline{a}) \subset {\rm
  Hilb}^p(\PP^{n})$ the stratum of standard determinantal schemes where
$f_{ij}$ are homogeneous polynomials of degrees $a_j-b_i$. Inside
$W_s(\underline{b};\underline{a})$ we have the open subset
$W(\underline{b};\underline{a})$ of determinantal schemes which are
generically a complete intersection. The elements are called {\it good}
determinantal schemes. Note that $W_s(\underline{b};\underline{a}) $ is
irreducible, and $W(\underline{b};\underline{a}) \ne \emptyset$ if we suppose
$a_{i-1}-b_i > 0$ for $i \ge 1$, see \eqref{WWs}.

In this paper we determine the dimension of a non-empty
$W(\underline{b};\underline{a})$ provided $a_{i-2}-b_i \ge 0$ for $i \ge 2$
and $n-c \ge 1$
(Theorem~\ref{Amodulethm3}, Corollary~\ref{Amodulethm4}). Indeed
 \begin{equation} \label{dimW}
\dim   W(\underline{b};\underline{a}) =  \lambda_c + K_3+K_4+...+K_c \, , 
 \end{equation}
 where $\lambda_c$ and $K_i$ are a large sum of binomials only involving $a_j$
 and $b_i$ (see Conjecture~\ref{conjnew1} and \eqref{lamda} for the definition
 of $\lambda_c$ and $K_i$). In terms of the Hilbert function, $H_M(-)$, of
 $M$, we may alternatively write \eqref{dimW} in the form $$\dim
   W(\underline{b};\underline{a})= \sum_{j=0}^{t+c-2} H_M(a_j)- \sum_{i=1}^{t}
   H_M(b_i) +1 $$(Remark~\ref{dimWba2}). Moreover we prove that the closure
 $\overline{ W(\underline{b};\underline{a})}$ is a generically smooth
 irreducible component of the Hilbert scheme $\Hi ^p(\PP^{n})$ provided $ _0\!
 \Ext^2_A(M,M)$, the degree zero part of the graded module $ \Ext^2_A(M,M)$,
 vanishes for a general $X=\Proj(A)$ of $ W(\underline{b};\underline{a})$
 (Theorem~\ref{Amodulethm5}, Corollary~\ref{Amodulethm6}).
Indeed
 \begin{equation*} \dim_{(X)} \Hi ^p(\PP^{n}) - \dim
   W(\underline{b};\underline{a}) \le \dim \ _0\! \Ext_A^2(M,M) \ ,
 \end{equation*}
 and $ \Hi ^p(\PP^{n})$ is smooth at $(X)$ if equality holds. We prove that
 $n-c \ge 2$ implies $ \Ext^2_A(M,M)=0$ (Corollary~\ref{Amodulecor4}), whence
 $\overline{ W(\underline{b};\underline{a})}$ is a generically smooth
 irreducible component of $\Hi ^p(\PP^{n})$ in the case $n-c \ge 2$ and
 $a_{i-\min (3,t)}\ge b_{i}$ for $\min (3,t)\le i \le t$. This proves
 Conjecture 4.2 of \cite{KM09}. Moreover our results hold for every $(X) \in
 W(\underline{b};\underline{a})$ provided a depth condition on the singular
 locus is fulfilled.\! A general $X$ of $ W(\underline{b};\underline{a})$
 satisfies the condition and we get the mentioned results.

 The most remarkable finding in this paper is perhaps the method. Indeed an
 embedded deformation problem for the determinantal scheme $X=\Proj(A)$,
 $A=R/{\rm ann}(M)$ is transfered to a deformation problem for the $R$-{\it
   module} $M$ where it is handled much more easily because every deformation
 of $M$ comes from deforming the matrix $\cA$. The latter is easy to see from
 the Buchsbaum-Rim complex. In fact it was in \cite{K09} we introduced the
 notion ``every deformation of $X$ comes from deforming $\cA$'' to better
 understand why $\overline{ W(\underline{b};\underline{a})}$ may fail to be an
 irreducible component. This led us to study deformations of $M$ because the
 corresponding property holds for $M$. Therefore the isomorphism between the
 graded deformation functors of $M$ and $R \to A$, which we prove under the
 assumption $ _0\! \Ext_A^i(M,M)=0$ for $i=1$ and $2$, is an important result
 (Theorem~\ref{compthm}). Note that the graded deformation functor of $R \to
 A$ is further isomorphic to the local Hilbert (scheme) functor of $X$
 provided $n-c \ge 1$. Since we also prove that $n-c \ge i \ge 1$ implies $
 \Ext^i_A(M,M)=0$ under mild assumptions (Theorems~\ref{Amodulethm} and
 \ref{Amodulethm2}), we get our rather algebraic method for studying a
 geometric object, the Hilbert scheme. Even the corresponding non-graded
 deformation functors of $M$ and $R \to A$ are isomorphic for $n-c \ge 2$
 (Remark~\ref{Amodulerem2}), which more than indicates that this method holds
 for local determinantal rings of dimension greater than $2$. Hence we expect
 applications to deformations of determinantal singularities, as well as to
 multigraded Hilbert schemes. We remark that while the vanishing of
 $\Ext_A^i(M,M)$ in Theorem~\ref{Amodulethm} is mainly known (at least for
 $i=1$, see Remark~\ref{SchSv}), 
 the surprise is Theorem~\ref{Amodulethm2} which reduces the depth assumption
 of Theorem~\ref{Amodulethm} by 1 in important cases. Note that the local
 deformation functors of $M$ as an $A$- as well as an $R$-module were
 thoroughly studied by R. Ile in \cite{I01}, \cite{IleTr} and in \cite{Ile}
 he studies the case of a determinantal hypersurface $X$ ($\cA$ a square
 matrix) without proving, to our knowledge,
 the mentioned  results (see Remark~\ref{rile}). 
 Ile and his paper \cite{Ile}, and the joint papers \cite{KMMNP}, \cite{KM},
 \cite{KM09} have, however, served as an inspiration for this work.

 We also get further interesting results, e.g. that arbitrary modules of
 maximal grade are unobstructed (earlier proved by Ile in \cite{I01}), and we
 show that the dimension of their natural deformation spaces is equal to the
 right hand side of \eqref{dimW} (Theorems~\ref{modulethm} and
 ~\ref{moduleCIthm}, cf.\! Remark~\ref{dimWba2}). Moreover we prove that the
 cohomology $H^i({\mathcal N}_X(v))$ of the normal sheaf of $X \subset
 \PP^{n}$ for a $X$ general in $ W(\underline{b};\underline{a})$ vanishes for
 $1 \le i \le \dim X-2$ and every $v$ (Theorem~\ref{normal}). Even the algebra
 cohomology groups ${\rm H}^{i}(R,A,A)$ of Andr\'{e}-Quillen vanish for $2
 \le i \le \min \{ \dim X-1, c \}$. This extends a result from T. Svanes'
 thesis \cite{Sv} proven there for so-called {\em generic} determinantal
 schemes in which the entries of $\cA$ are the indeterminates of $R$, see
 \cite{b-v}, Thm.\! 15.10 for details. Finally we remark that the assumption
 that $R$ is a polynomial ring can be weakened. Indeed all theorems and their
 proofs generalize at least to the case where $\Proj(R)$ is any arithmetically
 Cohen-Macaulay $k$-scheme (and smooth in Theorem~\ref{normal} (ii)), only
 replacing all $ \binom{v+n}{n}$ in \eqref{dimW} with $\dim R_{v}$.

 \vskip 1mm The method of this paper has the power of solving most of the
 deformation problems the author, together with coauthors (mostly Mir\'o-Roig
 at Barcelona) has considered in several papers (\cite{KMMNP}, \cite{KM},
 \cite{KM09}, \cite{K09}), mainly: \vskip 2mm (1) \
 Determine the dimension of $W(\underline{b};\underline{a})$ in terms of $a_j$
 and $b_i$ (see Conjecture~\ref{conjnew1}).

 (2) \ Is $\overline{ W(\underline{b};\underline{a})}$ a generically
 smooth irreducible component of $\Hi ^p(\PP^{n})$?

 \vskip 2mm The main method so far has been to delete columns of the matrix
 $\cA$, to get a ``flag'' of closed subschemes $X =X_c \subset X_{c-1} \subset
 ... \subset X_2 \subset {\PP^{n} }$ and to prove the
 results 
 by considering the smoothness of the Hilbert flag scheme of pairs and its
 natural projections into the Hilbert schemes. In fact in \cite{KM} we solved
 problem (1) in the cases $2\le c\le 5$ and $n-c \ge 1$ (assuming $char (k)=0$
 if $c=5$),
 and recently we almost solved (1) in the remaining cases under the assumption
 $ a_{t+3}> a_{t-2}$ \cite{KM09}. Concerning problem (2) we gave in \cite{KM}
 an affirmative answer in the range $2\le c \le 4$ and $n-c\ge 2$, (see
 \cite{elli} and \cite{KMMNP} for the cases $2 \le c \le 3$). We got further
 improvements in \cite{KM09} and conjectured a positive answer to problem (2)
 provided $n-c \ge 2$, but we were not able to solve all technical challenges
 which increased with the codimension. In this paper we fully prove the
 conjecture,
 as well as Conjecture~\ref{conjnew1} for $n-c \ge 1$, with the new approach
 which is much easier than the older one. For the case $n-c=0$ we remark that
 since every element of $ W(\underline{b};\underline{a})$ has the same Hilbert
 function, problem (2) becomes more natural provided we replace $\Hi
 ^p(\PP^n)$ with $\GradAlg(H)$, see \cite{K09} for details and the notations
 below.
 \vskip 2mm We thank R. Ile, R.M. Mir\'o-Roig, J.A. Christophersen, M. Boij,
 O.A. Laudal, Johannes Kleppe and U. Nagel for interesting discussions on
 different aspects of this topic.

\vskip 2mm

{\bf Notation:} In this work $R=k[x_0, \dots ,x_n]$ will be a polynomial ring
over an algebraically closed field, 
$\goth m= (x_0, \dots ,x_n)$ and $\deg x_i = 1$, unless explicitly making
other assumptions. 

We mainly keep the notations of \cite{KM} and \cite{K09}. If $X \subset Y$ are
closed subschemes of $\PP:=\PP^n:=\Proj(R)$, we denote by ${\mathcal
  I}_{X/Y}$ (resp. ${\mathcal N}_{X/Y}$) the ideal (resp. normal) sheaf of $X$
in $ Y$, and we usually suppress $Y$ when $Y=\PP^n$. By the codimension, ${\rm
  codim}_Y X$, of $X$ in $Y$ we simply mean $\dim Y -\dim X$, and we let
$\omega_X ={\mathcal E}xt^c_{{\mathcal O}_{\PP^n}} ({\mathcal O}_X,{\mathcal
  O}_{\PP^n})(-n-1)$ if $c= {\rm codim}_{\PP}
X$. 
When we write $X=\Proj(A)$ we take $A:=R/I_X$ and $K_A=\Ext^c_R (A,R)(-n-1)$
where $I_X=H^0_{*}({\mathcal I}_X)$ is the saturated homogeneous ideal of $X
\subset \PP^n$. We denote the group of morphisms between coherent
$\cO_X$-modules 
by $\Hom_{\cO_X}(\cF,\cG)$ while $\cH om_{\cO_X}(\cF,\cG)$ denotes the
corresponding sheaf. 
Moreover we set $\hom(\cF,\cG)=\dim_k\Hom(\cF,\cG)$ and we correspondingly use
small letters for the dimension, as a $k$-vector space, of similar groups. 

We denote the Hilbert scheme by $\Hi ^p(\PP^n)$, $p$ 
the Hilbert polynomial \cite{G}, and $(X) \in \Hi ^p(\PP^n)$ the point which
corresponds to $X\subset \PP^n$. Let $\GradAlg(H)$ be the representing object
of the functor which parametrizes flat families of graded quotients $A$ of $R$
of $\depth_{\goth m} A \ge 1$ and with Hilbert function $H$; $H(i)=\dim
A_i$ (\cite{K98}, \cite{K04}). We let $(A)$, or $(X)$ where $X= \Proj(A)$,
denote the point of $\GradAlg(H) $ which corresponds to $A$. Then $X$ (resp.\!
$A$) is {\it unobstructed} if $\Hi ^p(\PP^n)$ (resp.\! $\GradAlg(H) $) is smooth
at $(X)$. By \cite{elli},
\begin{equation} \label{Grad} \GradAlg(H) \simeq \Hi ^p(\PP^n) \ \ \ {\rm at}
  \ \ \ (X) \ \
\end{equation}
provided $\depth A:=\depth_{\goth m}A \geq 2$. This implies that if $\dim A
\ge 2$ and $A$ is Cohen-Macaulay (CM), then it is equivalent to consider
deformations of $X \hookrightarrow \PP^n$, or of $R \twoheadrightarrow A$, and
moreover that their tangent spaces $ _0\!\Hom (I_{X},A) \simeq \ H^0({\mathcal
  N}_X)$ are isomorphic where the lower index means the degree zero part of
the graded module $ \Hom (I_{X},A)$. We also deduce that if $X$
is generically a complete intersection, then $ _0\!\Ext^1_A(I_{X}/I_{X}^2,A)$
is an obstruction space of $ \Hi ^p(\PP^n)$ at $(X)$ (\cite{K04}, \S 1.1).
Finally we say that $X$ is {\it general} in some irreducible subset $ W
\subset \Hi ^p(\PP^n)$ if $(X)$ belongs to a sufficiently small open subset
$U$ of $W$ such that any $(X)$ in $U$ has all the openness properties that we
want to require.


\section{Background}

This section recalls basic results on standard and good determinantal schemes
needed in the sequel, see \cite{b-v}, \cite{eise}, \cite{BH} and \cite{KMNP}
for more details and \cite{e-n}, \cite{BR}, \cite{e-h} for background. Let
\begin{equation}\label{gradedmorfismo} \varphi:F=\bigoplus
  _{i=1}^tR(b_i)\longrightarrow G:=\bigoplus_{j=0}^{t+c-2}R(a_j)
\end{equation}
be a graded morphism of free $R$-modules and let
$\cA=(f_{ij})_{i=1,...t}^{j=0,...,t+c-2}$, $\deg f_{ij}=a_j-b_{i}$, be a
$t\times (t+c-1)$ homogeneous matrix which represents the dual
$\varphi^*:=\Hom_R(\varphi,R)$. Let $I_t(\cA)$  be
the ideal of $R$ generated by the maximal minors of $\cA$. In this paper we
suppose
$$ c\ge 2 \ , \ \ t\ge 2\ , \ \ b_1 \le ... \le b_t \ \ \ \ {\rm and} \ \ \ \
a_0 \le a_1\le ... \le a_{t+c-2}.$$ Recall that a codimension $c$ subscheme
$X\subset \PP^{n}$ is standard determinantal if $I_X=I_t(\cA)$ for some
homogeneous $t\times (t+c-1)$ matrix $\cA$ as above. Moreover $X\subset
\PP^{n}$ is a \emph{good determinantal} scheme if additionally, $I_{t-1}(\cA)$
defines a scheme of codimension greater or equal to $c+1$ in $ \PP^{n}$. Note
that if $X$ is standard determinantal and a generic complete intersection in $
\PP^{n}$, then $X$ is \emph{good determinantal}, and conversely \cite{KMNP},
Thm.\! 3.4. We say that $\cA$ is minimal if $f_{ij}=0$ for all $i,j$ with
$b_{i}=a_{j}$.

Let $W(\underline{b};\underline{a})$ (resp. $W_s(\underline{b};\underline{a}))$
be the stratum in $ {\rm Hilb}^p(\PP^{n})$ consisting of good (resp. standard)
determinantal schemes. 
By \cite{KM}, see the end of p.\! 2877, we get that the closures of these
strata in $ \Hi ^p(\PP^n)$ are equal and irreducible. Moreover since we will
not require $\cA$ to be minimal for $X=\Proj(R/I_{t}(\cA))$ to belong to
$W(\underline{b};\underline{a})$ or $W_s(\underline{b};\underline{a})$ in
their definitions (a slight correction to \cite{KMMNP} and \cite{KM}!), we
must reconsider Cor.\! 2.6 of \cite{KM}. Indeed we may use
its proof to see (cf. \cite{KM09} for details)
\begin{equation} \label{WWs} W(\underline{b};\underline{a}) \ne \emptyset \ \
  \Leftrightarrow \ \ W_s(\underline{b};\underline{a}) \ne \emptyset \ \
  \Leftrightarrow \ \ a_{i-1} \ge b_i \ \ {\rm for \ all} \ i \ {\rm and \ }
    a_{i-1} > b_i \ \ {\rm for \ some } \ i .
\end{equation}

Let $A=R/I_X$ (i.e. $X$) be standard determinantal and let $M:= \coker
(\varphi^*)$. Then one knows that the {\em Eagon-Northcott complex} yields the
following free resolution
\begin{equation}\label{EN}0 \ra \wedge^{t+c-1}G^* \otimes S_{c-1}(F)\otimes
  \wedge^tF\ra \wedge^{t+c-2} G ^*\otimes S _{c-2}(F)\otimes \wedge ^tF\ra
  \ldots  \end{equation}  
$$ \ra
\wedge^{t}G^* \otimes S_{0}(F)\otimes \wedge^tF\ra R \ra A \ra 0 $$ of $A$ and
that the {\em Buchsbaum-Rim complex} yields a free resolution of $M$;
\begin{equation}\label{BR}
0 \ra \wedge^{t+c-1}G^* \otimes S_{c-2}(F)\otimes \wedge^tF\ra
\wedge^{t+c-2} G ^*\otimes S _{c-3}(F)\otimes \wedge ^tF\ra \ldots 
\end{equation}
 $$ \ra \wedge^{t+1}G^* \otimes S_{0}(F)\otimes \wedge^tF \ra G^*\ra \ F^* \ra
 M \ra 0 ,$$ (the resolutions are minimal if $\cA$ is minimal), see for
 instance \cite{b-v}; Thm.\! 2.20 and \cite{eise}; Cor.\! A2.12 and Cor.\!
 A2.13. Note that \eqref{EN} shows that A is Cohen-Macaulay. 

 Let ${\cB}$ be the matrix obtained by deleting the last column of ${\cA}$ and
 let $B$ be the $k$-algebra given by the maximal minors of ${\cB}$. Let
 $Y=\Proj(B)$. The transpose of ${\cB}$ induces a map $ \phi:F=\oplus
 _{i=1}^tR(b_i)\rightarrow G':=\oplus _{j=0}^{t+c-3}R(a_j)$. Let $M_{\cB}$ be
 the cokernel of $\phi^*=\Hom_R(\phi,R)$ and let $M_{\cA}=M$ and $c >
 2$.
In this situation we recall that there is an exact sequence
\begin{equation}\label{Mi}
0\longrightarrow B \longrightarrow M_{\cB}(a_{t+c-2})
\longrightarrow M_{\cA}(a_{t+c-2}) \longrightarrow 0
\end{equation}
in which $B \longrightarrow M_{\cB}(a_{t+c-2})$ is a regular section given by
the last column of ${\cA}$. Moreover,
\begin{equation}\label{Di}
0\longrightarrow
M_{\cB}(a_{t+c-2})^* :=\Hom_{B}(M_{\cB}(a_{t+c-2}),B)\longrightarrow B
\longrightarrow A \longrightarrow 0
\end{equation}
is exact by \cite{KMNP} or \cite{KMMNP}, (3.1), i.e. we may put
$I_{X/Y}:=M_{\cB}(a_{t+c-2})^*$. Due to \eqref{BR}, $M$ is a maximal
Cohen-Macaulay $A$-module ($\depth M = \dim A$), and so is $I_{X/Y}$ by
(\ref{Di}). By \cite{eise} we have $K_{A}(n+1)\cong S_{c-1}M_{\cA}(\ell_c)$
and hence $K_{B}(n+1)\cong S_{c-2}M_{\cB}(\ell_{c-1})$ where
\begin{equation}\label{ell}
 \ell_i :=\sum_{j=0}^{t+i-2}a_j-\sum_{k=1}^tb_k \ \ {\rm for} \ \ 2 \le i \le c.
\end{equation}
Recall that $\tilde{M}$ is locally free of rank one precisely on
$X-V(I_{t-1}(\cA))$ (\cite{BH}, Lem.\! 1.4.8) and that $X\hookrightarrow
\PP^{n}$ is a local complete intersection (l.c.i.) by e.g. \cite{ulr}, Lem.\!
1.8 provided we restrict to $X-V(I_{t-1}(\cA))$. 
By \eqref{Di} it follows that $X\hookrightarrow Y$ and $Y\hookrightarrow
\PP^{n}$ are l.c.i.'s outside $V(I_{t-1}(\cB))$. Note that $V(I_{t-1}(\cB))
\subset V(I_{t}(\cA))=X$.

\begin{remark} \label{dep} 
  Put $X_c:=X$ and $X_{c-1}:=Y$, let $c > 2$ and let $\alpha$ be a positive
  integer. If $X$ is general in $W(\underline{b};\underline{a})$ and
  $a_{i-\min (\alpha ,t)}-b_i\ge 0$ for $\min (\alpha
  ,t)\le i\le t$, then
\begin{equation}\label{a-b} \codim_{X_j}Sing(X_j)\ge \min\{2\alpha
  -1, j+2\}  \ \ {\rm for \  \ } j=c-1 {\rm \ and \ } c \ . \end{equation}
This follows from Rem.\! 2.7 of \cite{KM} (i.e., from \cite{chang}).
In particular if $\alpha \ge 3$, we get that the closed embeddings $Y
\hookrightarrow \PP^{n}$ and $X \hookrightarrow Y$ are local complete
intersections outside some set $Z$ of codimension at least $ \min(4,c)$.
Indeed we may take $Z=V(I_{t-1}(\cB))$.
\end{remark}

Moreover we recall the following useful general fact that if $L$ and $N$ are
finitely generated $A$-modules such that $\depth_{I(Z)}L\ge r+1$ and
$\tilde{N}$ is locally free on $U:=X-Z$, then the natural map
\begin{equation} \label{NM}
\Ext^{i}_A(N,L)\longrightarrow
H_{*}^{i}(U,{\cH}om_{{\cO}_X}(\tilde{N},\tilde{L}))
\end{equation}
is an isomorphism, (resp. an injection) for $i<r$ (resp. $i=r$), and
$H_{*}^{i}(U,{\cH}om_{{\cO}_X}(\tilde{N},\tilde{L}))\simeq
H^{i+1}_{I(Z)}(\Hom_A(N,L))$ for $i>0$, cf.\! \cite{SGA2}, exp.\! VI. Note
that we interpret $I(Z)$ as $\goth m$ if $Z= \emptyset$.
\vskip 2mm
In \cite{KM} 
we conjectured the dimension of  $
W(\underline{b};\underline{a})$ in terms of the invariant {\small
\begin{equation} \label{lamda} \lambda_c:= \sum_{i,j}
    \binom{a_i-b_j+n}{n} + \sum_{i,j} \binom{b_j-a_i+n}{n} - \sum _{i,j}
    \binom{a_i-a_j+n}{n}- \sum _{i,j} \binom{b_i-b_j+n}{n} + 1. 
  \end{equation} }
Here the indices belonging to $a_j$ (resp. $b_i$) range over $0\le j \le
t+c-2$ (resp. $1\le i \le t$) and we let $ \binom{a}{n}= 0$ if $a$ is a
negative integer. Since \cite{K09} discovered that the scheme
of $c+1$ general points in $\PP^{c}$ given by the vanishing all $2 \times 2$
minors of a general $2 \times (c+1)$ 
matrix of linear entries was a counterexample to
Conjecture 6.1 (and to the special case given in Conjecture 6.2) of \cite{KM},
we  slightly changed  Conjecture 6.1 in \cite{KM09} to
\begin{conjecture} \label{conjnew1} Given integers $a_0\le a_1\le ... \le
  a_{t+c-2}$ and $b_1\le ...\le b_t$, let  $h_{i-3}:= 2a_{t+i-2}-\ell_i +n$,
  for $i=3,4,...,c$ and assume $a_{i-\min ([c/2]+1,t)} \ge
  b_{i}$ provided $n>c$ and $a_{i-\min ([c/2]+1,t)} >
  b_{i}$ provided $n=c$
  for $\min ([c/2]+1,t)\le i \le t$. Except for the family
  $W({0,0};{1,1,...,1})$ of \ {\rm zero dimensional} schemes above we have,
  for $W(\underline{b};\underline{a}) \ne \emptyset$, that
 \[
 \dim W(\underline{b};\underline{a}) = \lambda_c+ K_3 + K_4+...+K_c \ , \]
 where $K_3=\binom{h_0}{n}$ and $K_4= \sum_{j=0}^{t+1} \binom{h_1+a_j}{n}-
 \sum_{i=1}^{t} \binom{h_1+b_i}{n}$ and in general \[ K_{i+3}= \sum _{r+s=i
   \atop r , s \ge 0} \sum _{0\le i_1< ...< i_{r}\le t+i \atop 1\le
   j_1\le...\le j_s \le t } (-1)^{i-r} \binom{h_i+a_{i_1}+\cdots
   +a_{i_r}+b_{j_1}+\cdots +b_{j_s} }{n} \ {\rm for} \ 0 \le i \le c-3. \]
\end{conjecture}
In \cite{KM}, Thm.\! 3.5 we proved that the right hand side in the formula for
$\dim W(\underline{b};\underline{a})$ in the Conjecture is always an upper
bound for $\dim W(\underline{b};\underline{a})$, and moreover,
that the Conjecture hold in the range
\begin{equation} \label{2c5}
2 \le c \le 5 \  \
  {\rm and} \ \ n-c \ge 1 \ {\rm (\, supposing \ } char k = 0 {\rm \ if} \ c =
    5\, )\, .
  \end{equation}
  Indeed this is mainly \cite{KM}, Thm.\! 4.5, Cor.\! 4.7, Cor.\! 4.10, Cor.\!
  4.14 and \cite{elli} ($c=2$) and \cite{KMMNP} ($c=3$). Moreover we have by
  \cite{KM09} (valid also for $n=c$ without assuming $char k = 0$):

\begin{theorem} \label{dim0new} Assume $a_{0} > b_{t}$. Then
  Conjecture~\ref{conjnew1} holds provided $c>5$ (resp. $2 \le c \le 5$) and $
  \ a_{t+3} > a_{t-2}$ (resp.$ \ a_{t+c-2} > a_{t-2} $).
\end{theorem}

In \cite{KM09} we stated a Conjecture related to the problem (2) 
of the Introduction:
\begin{conjecture} \label{conjcomp} Given integers $a_0\le a_1\le ... \le
  a_{t+c-2}$ and $b_1\le ...\le b_t$, we suppose $n-c\ge 2$, $c \ge 5$ and
   $ \ a_{0} > b_{t}$. 
Then $\overline{
    W(\underline{b};\underline{a})}$ is a generically smooth irreducible
  component of the Hilbert scheme $\Hi ^p(\PP^{n})$.
\end{conjecture}
By \cite{KM09}, Cor.\! 3.8 and Thm.\! 3.4, Conjecture~\ref{conjcomp} holds
provided $ \ a_{t+3} > a_{t-1} + a_{t}-b_1$ or more generally, if a certain
collection of $\Ext^1$-groups vanishes. Note that the conclusion of
Conjecture~\ref{conjcomp} holds if $n-c\ge 2$ and $2 \le c \le 4$ by
\cite{elli}, \cite{KMMNP} and \cite{KM}. 

\vspace{1.5mm} As in \cite{K09} we briefly say ``$T$ a local ring'' (resp.
``$T$ artinian'') for a local $k$-algebra $(T,{\mathfrak m}_T)$ essentially of
finite type over $k=T/{\mathfrak m}_T$ (resp. such that ${\mathfrak m}_T^r=0$
for some integer $r$). The local deformation functors of this paper will be
defined over the category $\underline {\ell}$ of artinian $k$-algebras.
Moreover we say ``$T \twoheadrightarrow S$ is a small in $\underline {\ell}$
'' provided there is a surjection $(T,{\mathfrak m}_T) \rightarrow (S,
{\mathfrak m}_S)$ of artinian $k$-algebras whose kernel ${\mathfrak a}$
satisfies ${\mathfrak a} \cdot {\mathfrak m}_T=0$.

If $T$ is a local ring, we denote by $\cA_T=(f_{ij,T})$ a matrix of
homogeneous polynomials belonging to the graded polynomial algebra $R_T:=R
\otimes_k T$, satisfying $f_{ij,T} \otimes_T k=f_{ij}$ and $\deg f_{ij,T}
=a_j-b_i$. 
Note that all elements from $T$ are considered to be of degree zero. For short
we say $\cA_T$ {\em lifts $\cA$ to} $T$. The matrix $\cA_T$ induces a morphism
\begin{equation} \label{al} \varphi_T: F_T:=\oplus _{i=1}^tR_T(b_i)\rightarrow
  G_T:=\oplus_{j=0}^{t+c-2}R_T(a_j) \ .
\end{equation}
\begin{lemma} \label{defalpha} If $X=\Proj(A)$, $A=R/I_t(\cA)$, is a standard
  determinantal scheme, then $A_T:=R_T/I_t(\cA_T)$ and $M_T:= \coker
  \varphi_T^*$ are (flat) graded deformations of $A$ and $M$ respectively for
  every choice of $\cA_T$ as above. In particular $X_T=\Proj(A_T) \subset
  \PP^n_T:=\Proj(R_T) $ is a deformation of $X \subset \PP^n$ to $T$ with
  constant Hilbert function.
\end{lemma}
\begin{proof} [Proof {\rm(}\cite{K09}, {\rm cf.} \cite{Sha}, {\rm Rem. to
    Prop.\! 1)}] 
  The Eagon-Northcott and Buchsbaum-Rim complexes are functorial in the sense
  that, over $R_T$, all free modules and all morphisms in these complexes
are determined by $\cA_T$. Since these complexes become free resolutions of
$A$ and $M$ respectively when we tensor with $k$ over $T$, it follows that
$A_T$ and $M_T$ are flat over $T$ and satisfy $A_T \otimes_T k = A$ and $M_T
\otimes_T k = M$.
\end{proof}
\begin{definition} \label{everydef} 
  We say ``every deformation of $X$ 
  comes from deforming $\cA$'' if for every local ring $T$ and every graded
  deformation $R_T \to A_T$ of $R \to A$ to $T$, then $A_T$ is of the form
  $A_T=R_T/I_t(\cA_T)$ for some $\cA_T$ as above. Note that by \eqref{Grad} we
  can in this definition replace ``graded deformations of $ R \to A$'' by
  ``deformations of $X \hookrightarrow \PP^{n}$'' provided $\dim X \ge 1$.
\end{definition}
\begin{lemma} \label{unobst} Let $X=\Proj(A)$ be a standard determinantal
  scheme, $(X) \in W(\underline{b};\underline{a})$. If every deformation of
  $X$ comes from deforming $\cA$, then $A$ is unobstructed. Moreover if $n-c
  \ge 1$ then $X$ is unobstructed and $\overline{
    W(\underline{b};\underline{a})}$ is an irreducible component of $\Hi ^p
  (\PP^{n}).$
\end{lemma}
\begin{proof} Let $T \twoheadrightarrow S$ be a small in $\underline {\ell}$
  and let $A_S$ be a deformation of $A$ to $S$. By assumption,
  $A_S=R_S/I_t(\cA_S)$ for some matrix $\cA_S$. We can lift each $f_{ij,S}$ to
  a polynomial $f_{ij,T}$ with coefficients in $T$ such that $f_{ij,T}
  \otimes_T S=f_{ij,S}$. By Lemma~\ref{defalpha}, $A_T:=R_T/I_t(\cA_T)$ is
  flat over $T$. Since $A_T \otimes_T S = A_S$ we get the unobstructedness of
  $A$, as well as the unobstructedness of $X$ in the case $\dim X \ge 1$ by
  \eqref{Grad}. For the remaining part of the proof, see \cite{K09}.
\end{proof}

\begin{remark} By these lemmas we get $T$-flat determinantal schemes by just
  parameterizing the polynomials of
  $\cA$ over a local ring $T$, see Rem.\! 4.5 of \cite{K09} and Laksov's
  papers \cite{dan2}, \cite{dan} for somewhat similar results for more general
  determinantal schemes. 
%
%
 %
%
%
\end{remark}

\section{deformations of $R$-modules of maximal grade}

Let $M$ be a finitely generated (torsion) $R$-module with presentation matrix
$\cA$, i.e. $M=\coker (\varphi^*)$ with $\varphi$ as in
\eqref{gradedmorfismo}. Since the grade of $M$ over $R$ is the grade, or
codimension, of the annihilator $I:= {\rm ann}(M)$ of $M$, and since the
radicals of $ I$ and $I_t(\cA)$ are the same, we get that $M$ has maximal
grade if and only if $A:=R/I_t(\cA)$ is standard determinantal. In this case $
I =I_t(\cA)$, see \cite{eis} for details. If $M \cong R/I(-b_1)$ is cyclic
($t=1$), we remark that a module of maximal grade is a complete intersection.
The main results of this section is variations of the following

\begin{theorem} \label{modulethm} Let $M$ be a finitely generated graded
  module over $R$ of maximal grade. Then $M$ is unobstructed. Moreover if
  $A:=R/{\rm ann}(M)$ is generically a complete intersection, then
 $$\dim \ _0\! \Ext_R^1(M,M) = \lambda_c + K_3+K_4+...+K_c \ \ \ {\rm and}
$$ \\[-8mm]
 $$\depth \Ext_R^1(M,M) \ge \dim A -1 \ .$$ 
\end{theorem}

\begin{remark} \label{modulerem} By deformation theory $\ _0\! \Ext_R^1(M,M)$
  (resp. $ \ _0\! \Ext_R^2(M,M)$) is the tangent (resp. the natural
  obstruction) space of the local deformation functor, ${\rm Def}_{M/R}$, of
  $M$ as a graded $R$-module (e.g. \cite{siq}). Since 
  $c \ge 2$, $ \ _0\! \Ext_R^2(M,M)$ is in many cases non-vanishing.
\end{remark}
\begin{remark} \label{modulerem0} Note that the assumption on $A$ in
  Theorem~\ref{modulethm} is equivalent to assuming $A$ good determinantal. By
  \eqref{WWs} good determinantal schemes exist if standard
  determinantal schemes exist. Hence if we take the polynomials $f_{ij}$ of
  degrees $a_j-b_i$ in a presentation matrix $(f_{ij})$ of $M$ general enough,
  then the assumption on $A$ in Theorem~\ref{modulethm} is satisfied.
%
\end{remark}
\begin{remark} \label{ile} While distributing a preliminar version of a paper
  partially containing Theorem~\ref{modulethm} to specialists in deformations
  of modules, we learned that the unobstructedness part of
  Theorem~\ref{modulethm} (and hence of Theorem~\ref{moduleCIthm}) was proved
  in R. Ile's PhD thesis, cf.\! \cite{I01}, ch.\! 6.
\end{remark}
\begin{proof} Let $T \twoheadrightarrow S$ be a small in $\underline {\ell}$
  and let $M_S$ be any graded deformation of $M$ to the artinian ring $S$. Let
  $\cA =(f_{ij})$ be a homogeneous matrix which represents $\varphi^*$. Since
  $G^* \stackrel {\varphi^*} { \longrightarrow} F^* \to M \to 0$ is exact (cf.
  \eqref{gradedmorfismo}), we have $M_S=\coker (\varphi_S^*)$ where
  $\varphi_S^*$ corresponds to some matrix $\cA_S=(f_{ij,S})$, as in
  \eqref{al}. Since $T \rightarrow S$ is surjective, we can lift each
  $f_{ij,S}$ to a polynomial $f_{ij,T}$ with coefficients in $T$ such that
  $f_{ij,T} \otimes_T S=f_{ij,S}$. By Lemma~\ref{defalpha}, $M_T:=\coker
  (\varphi_T^*)$ is flat over $T$ and since $M_T \otimes_T S = M_S$ it follows
  that $M$ is unobstructed.

 To see the dimension formula we {\it claim} that there is an exact sequence
\begin{equation} \label{homext}
  0 \to  \ _0\! \Hom_R(M,M) \to\ _0\! \Hom_R(F^*,M) \to \ _0\!
  \Hom_R(G^*,M) \to \ _0\! \Ext_R^1(M,M) \to 0.
\end{equation}
Indeed look at the map $d_1: \wedge^{t+1}G^* \otimes S_{0}(F)\otimes \wedge^tF
\to G^*$ appearing in the Buchsbaum-Rim complex \eqref{BR} and recall that the
image of the corresponding map $\wedge^{t}G^* \otimes S_{0}(F)\otimes
\wedge^tF \to R$ of the Eagon-Northcott complex \eqref{EN} is the ideal $I=
{\rm ann}(M)$ generated by the maximal minors. It follows that $\im d_1
\subset I \cdot G^*$ and hence that the induced map $\ _0\! \Hom_R(d_1,M)=0$.
So if we apply $\ _0\! \Hom_R(-,M)$ to \eqref{BR}, we get \eqref{homext} by
the definition of $\ _0\! \Ext_R^i(M,M)$.

Let $E=\coker \varphi$. Then we have an exact sequence $$ 0 \ra E^* \ra G^*
\stackrel {\varphi^*} { \longrightarrow} F^* \ra M \ra 0 \ ,$$ to which we
apply the exact functors $ _0\! \Hom_R(F^*,-)$ and $ _0\! \Hom_R(G^*,-)$.
We get
\[ \ _0\! \hom_R(G^*,M) - \ _0\! \hom_R(F^*,M) \ = \lambda_c - 1 + \ _0\!
\hom_R(G^*,E^*) - \ _0\! \hom_R(F^*,E^*) \]
by using the definition \eqref{lamda} of $ \lambda_c$. Note that $\ _0\!
\hom(M,M)=1$ by \cite{KM}, Lem.\! 3.2 since $A$ is good determinantal by
assumption. Hence we get the dimension formula of Theorem~\ref{modulethm} from
\eqref{homext} provided we can prove $$ \ _0\! \hom_R(G^*,E^*) - \ _0\!
\hom_R(F^*,E^*) = K_3+...+K_c.$$ By \cite{KM}, Prop.\! 3.12 we have
$1+K_3+...+K_c= \ _0\! \hom_R(E,E)$ and by the proof of the same proposition
we find $ _0\! \hom_R(E,E)=1+ \ _0\! \hom_R(E,G) - \ _0\! \hom_R(E,F)$ and
whence we get the dimension formula.

Now we consider the depth of $ \Ext_R^1(M,M)$. Firstly observe that it is
straightforward to see $\depth \Ext_R^1(M,M) \ge \dim A -2$. Indeed using $
\Hom_R(M,M) \simeq A$ (\cite{KM}, Lem.\! 3.2) and skipping the lower index $0$
in \eqref{homext}, we get that all three $\Hom$-modules in \eqref{homext} are
maximal CM $A$-modules. Then if $D:= \coker(\Hom_R(M,M) \to \Hom_R(F^*,M))$ we
easily conclude since $\depth D \ge \dim A -1$ and $\depth \Ext_R^1(M,M) \ge
\depth D -1$ by \cite{eise}, Cor.\! 18.6.

Looking more carefully at the argument, we can show $\depth \Ext_R^1(M,M)
\ge \dim A -1$. Indeed it suffices to prove $\depth D = \dim A$. 
To see it we use the resolution of $A$ in \eqref{EN} and the resolution of $M
\otimes F$ deduced from \eqref{BR}. Let $\{f_1, f_2,...,f_t \}$ be the
standard basis of $F$ and $\{y_1, y_2,...\}$ the standard basis of $G^*$. The
leftmost free modules in these resolutions are $ \wedge^{t+c-1}G^* \otimes
S_{c-1}(F)\otimes \wedge^tF $ and $ \wedge^{t+c-1}G^* \otimes
S_{c-2}(F)\otimes \wedge^tF \otimes F $ respectively. We may consider the
former as an $R$-submodule of the latter through the
map 
$\tau_{c-1}$ where $\tau_k =id \otimes \tau'_{k}$ and $\tau'_k: S_{k}(F) \to
S_{k-1}(F)\otimes F$ is induced by sending a symmetric tensor $(f_{i_1}
\otimes_s...\otimes_s f_{i_k}) \in S_{k}(F)$ onto the ``reduced sum'' of
$\sum_{j=1}^k (f_{i_1} \otimes_s...\otimes_s \hat{f_{i_j}}
\otimes_s...\otimes_s f_{i_k})\otimes f_{i_j} \in S_{k-1}(F)\otimes F$. Here
``reduced'' means sending e.g. $(f_1 \otimes_s f_1 \otimes_s f_1 \otimes_s f_2
\otimes_s f_2 \otimes_s f_3) \in S_{6}(F)$ onto {\tiny
$$a(f_1 \otimes_s f_1 \otimes_s
f_2 \otimes_s f_2 \otimes_s f_3)\otimes f_1 +b(f_1 \otimes_s f_1 \otimes_s f_1
\otimes_s f_2 \otimes_s f_3) \otimes f_2 + (f_1 \otimes_s f_1 \otimes_s f_1
\otimes_s f_2 \otimes_s f_2) \otimes f_3$$} \\[-5mm] in $S_{5}(F) \otimes F$
with $(a,b)=(1,1)$ (and not
$(a,b)=(3,2)$!). 
Then $\tau_1 = id$ and letting $\tau_{-1}:R \to F^* \otimes F$ be the obvious
map and $\tau_0 : \wedge^{t}G^* \otimes S_{0}(F)\otimes \wedge^tF \to G^*
\otimes F $ the map induced by sending $(y_{i_1} \wedge...\wedge y_{i_t}) \in
\wedge^{t}G^*$ onto $\sum_{j=1}^t (-1)^j(\varphi^*(y_{i_1}) \wedge...\wedge
\hat{\varphi^*}(y_{i_j}) \wedge...\wedge \varphi^*(y_{i_t})) \otimes y_{i_j}
\in \wedge^{t-1}F^* \otimes G^*$ followed by the natural map $ \wedge^{t-1}F^*
\otimes G^* \to F \otimes \wedge^tF^* \otimes G^*$, one may check that the
collection of maps $\{ \tau_i \}_{i \ge -1}$ is actually a map between the
free resolutions of $A$ and $M \otimes F$. The explicit description in
\cite{Kir} of the differentials in the resolutions \eqref{EN} and \eqref{BR}
may be helpful in checking that the diagrams between the resolutions commute.
Now using the well known mapping cone construction we find a free resolution
of $D$ and since $\tau_{c-1} \otimes_R R/\goth m$ is injective, the leftmost
term $ \wedge^{t+c-1}G^* \otimes S_{c-1}(F)\otimes \wedge^tF $ becomes
redundant. The minimal $R$-free resolution of $D$ has therefore the same
length as the minimal $R$-free resolution of $A$, i.e. $D$ is maximally CM and
we are done. 
\end{proof}
\begin{remark} \label{modulerem1} We see from the proof, or
  Lemma~\ref{defalpha}, that if we arbitrarily lift the polynomials in a
  presentation matrix of $M$ to polynomials with coefficients in $T$, we get
  that $M_T:=\coker (\varphi_T^*)$ is {\it flat} over $T$. This is not true in
  general, but for modules of maximal grade it is because the Buchsbaum-Rim
  complex provides us with a resolution of $M_T$.
\end{remark}
\begin{remark} \label{isso} Let $E= \coker \varphi$, $ \varphi= F \to G$ \ cf.
  \eqref{gradedmorfismo} and suppose $R/I_t(\cA)$ is good
  determinantal. It is stated in \cite{KM}, Rem.\! 3.14 that $\dim \ _0\!
  \Ext_R^1(E,E) = \lambda_c + K_3+K_4+...+K_c$. Indeed one may use the proof
  above to see $ \ _0\! \Ext_R^i(E,E) \cong \ _0\! \Ext_R^i(M,M)$ for $i=1$
  while this is not true in general for $i \ne 1$.
\end{remark}
\begin{remark} \label{modulerem2} The theorem admits a {\it vast
    generalization} since the assumption that $R$ is a polynomial ring is not
  necessary. Indeed if $R$ is any commutative graded (resp. local) $k$-algebra,
  then a module of maximal grade is unobstructed and the exact sequence
  \eqref{homext} (resp. where the lower index 0 is removed) holds. In fact all
  we need for these parts in the proof is the existence and exactness of the
  Buchsbaum-Rim complex, which hold under almost no assumption on $R$ (cf.
  \cite{eise}, Appendix 2).
\end{remark} 

We will give the details in the graded case of what we claimed in
Remark~\ref{modulerem2}. This means that we will 
generalize to {\it arbitrary} modules of {\it maximal grade} the following
well-known fact for cyclic modules, that a complete intersection $R/I$ is
unobstructed and that $ _0\! \hom_R(I,R/I) = \sum_{j=1}^q \dim (R/I)_{(a_j)}$
where $\oplus_{j=1}^q R(-a_j) \to I$ is a minimal surjection.

For the remaining part of this section we let $R= \oplus_{v \ge 0} R_v$ be
{\it any} graded $k$-algebra ($k=R_0$ a not necessarily algebraically closed
field), generated by finitely many elements from $\goth m := \oplus_{v \ge 1}
R_v$. 
Let 
\begin{equation} \label{minpres}
G^*:= \sum_{j=1}^q R(-a_j) \stackrel {\varphi^*} { \longrightarrow} F^*:=
\sum_{j=1}^p R(-b_i)  \stackrel {\pi} { \longrightarrow} M \to 0
\end{equation}
be a minimal presentation of $M$ and suppose $M$ is of {\it maximal grade}.
Let $N= \ker \pi $. It is known that the tangent space of the graded
deformation functor ${\rm Def}_M(F^*)$ 
which deforms the surjection $ \pi: F^* \to M $ to artinian $k$-algebras from 
$\underline {\ell}$, using trivial deformations of $ F^*$, is isomorphic to $
_0\! \Hom_R(N,M)$ and that $ _0\! \Ext^1_R(N,M)$ contains all the obstructions
of the graded deformations (we may deduce it from \cite{L}, Thm.\! 4.1.14 and
Lem.\! 3.1.7, but \cite{G}, Prop.\! 5.1 and Cor. 5.2 and 5.3 is the classical
reference since we here deal with the local deformation functor, adapted to
graded deformations, of Grothendieck's Quot scheme). If we apply $ _0\!
\Hom_R(-,M)$ to $0 \to N \to F^* \to M \to 0$, we get the exact sequence
\begin{equation} \label{homext2}
  0 \to  \ _0\! \Hom_R(M,M) \to\ _0\! \Hom_R(F^*,M) \to \ _0\!
  \Hom_R(N,M) \to \ _0\! \Ext_R^1(M,M) \to 0 \ ,
\end{equation}
and $ _0\! \Ext^1_R(N,M) \simeq \ _0\! \Ext_R^2(M,M)$. We notice that the
arguments in the proof of Theorem~\ref{modulethm} which led to $\ _0\!
\Hom_R(d_1,M)=0$, where now $d_1: \wedge^{p+1}G^* \otimes S_{0}(F)\otimes
\wedge^pF \to G^*$, carry over to the general situation we are considering
since they 
relied on how the maps in the Buchsbaum-Rim complex were defined. Hence we get
the exact sequence \eqref{homext}, and comparing with \eqref{homext2}, we get
that the tangent space of ${\rm Def}_M(F^*)$ is $$\ _0\! \Hom_R(N,M) \simeq \
_0\! \Hom_R(G^*,M) \simeq \bigoplus_{j=1}^{q} M_{(a_j)}. $$ Since the map
$d_1$ is defined in terms of $ {\varphi^*} $, the unobstructedness argument
for $M$ in the proof of Theorem~\ref{modulethm} and the flatness argument of
Lemma~\ref{defalpha} both carry over the general case. Note that also the
object $ \pi:F^* \to M $ is unobstructed, i.e. ${\rm Def}_M(F^*)$ is formally
smooth (\cite{S}, \cite{L}) because $ \pi $ is easily deformed once $M$ is
deformed. Also the proof of the length of an $R$-free resolution of $
\Ext_R^1(M,M)$ holds and we have 
\begin{theorem} \label{moduleCIthm} Let $M$ be a finitely generated graded
  module (as in \eqref{minpres}) of maximal grade over a finitely generated
  graded $k$-algebra $R$ where $R_0=k$ is an arbitrary field. Let $N:= \ker (
  F^* \to M )$. Then $ \ \pd_R \Ext_R^1(M,M) \le c+1$. Moreover $M$ is
  unobstructed. Indeed ${\rm Def}_M(F^*)$ is formally smooth and the dimension
  of the tangent space of \ ${\rm Def}_M(F^*)$ is
 $$\dim \, _0\! \Hom_R(N,M) = \sum_{j=1}^{q} \dim M_{(a_j)}. $$
\end{theorem}

\begin{remark} \label{dimWba2} Under the assumptions of
  Theorem~\ref{moduleCIthm}, we see from the proof that 
 $$\dim \ _0\! \Ext_R^1(M,M)- \dim \, _0\! \Hom_R(M,M) = \sum_{j=1}^{q} \dim
 M_{(a_j)}- \sum_{i=1}^{p} \dim M_{(b_i)}.$$ Now suppose $R$ is any graded
 Cohen-Macaulay quotient of a polynomial ring $k[x_0, \dots ,x_n]$ with the
 standard grading where $k$ is any field. This will be the natural setting,
 having algebraic geometry in mind, to which we can generalize the theorems of
 this paper (sometimes assuming $k=\overline k$ to be algebraically closed).
 Slightly generalizing \cite{KM}, Lem.\! 3.2, we get that $ \Hom_R(M,M) \simeq
 A$ if $ \depth_{I_{t-1}(\cA)A}A \ge 1$ (cf. Remark~\ref{Amodulerem}). Hence $
 \, _0\! \hom_R(M,M)=1$ and the formula above gives an alternative to the
 formula of Theorem~\ref{modulethm} for computing $ _0\! \ext_R^1(M,M)$. In
 this general setting one may see that also the formula of
 Theorem~\ref{modulethm} holds provided we redefine $\lambda_c$ and $K_i$
 appearing in \eqref{lamda} and Conjecture~\ref{conjnew1} in the obvious way,
 namely by replacing all $ \binom{v+n}{n}$ with $\dim R_{v}$. Indeed this
 follows from the proof of Theorem~\ref{modulethm} since the part we use from
 \cite{KM} (Prop.\! 3.12) also generalize to this setting.
\end{remark}
\vskip 2mm

\section{the rigidity of modules of maximal grade}

In this section we consider a module $M$ of maximal grade as a graded module
over $A=R/I$ where $I={\rm ann}(M)$. 
Recall that $I=I_t(\cA)$ where ${\cA}=(f_{ij})_{i=1,...t}^{j=0,...,t+c-2}$ is
a $t\times (t+c-1)$ homogeneous presentation matrix of $M$, cf.
\eqref{gradedmorfismo}, in which case we put $M=M_{\cA}$. A main result in
this section is the rigidity of $M$ as an $A$-module (i.e. $\
\Ext_A^1(M,M)=0$) provided $X:=\Proj(A)$ is smooth of dimension greater or
equal to 1. Furthermore if $ \dim X \ge 2$ we also show $\ \Ext_A^2(M,M)=0$.
More generally we have the following results.
\begin{theorem} \label{Amodulethm} Let $M$ be a finitely generated graded
  $R$-module of maximal grade and let $A:=R/{\rm ann}(M)$. Let $j \ge 1$ be an
  integer and suppose $\depth_{I_{t-1}(\cA)A}A \ge j+1$. 
  Then $\Hom_A(M,M) \simeq A$ and,
 $$  \Ext_A^i(M,M) = 0 \ \ \ {\rm for} \ \ 1 \le i \le j-1.$$
  \end{theorem}

  \begin{remark} Let $X=\Proj(A)$,  $J:=I_{t-1}(\cA)$ and recall that
    $\depth_{JA}A = \dim A - \dim A/JA$ and $V(JA) \subset Sing(X)$. We may
    therefore take $j$ as $j= \codim_{X}Sing(X)-1= \dim X -\dim Sing(X)-1$ in
    Theorem~\ref{Amodulethm}, interpreting $\dim Sing(X)$ as $-1$ 
    if $ Sing(X)= \emptyset $.
 \end{remark}
 \begin{proof} 
   Since $\tilde M$ is locally free of rank one over
   $U:=\Proj(A)-V(I_{t-1}(\cA)A)$ (see the text before Remark~\ref{dep}), we
   can use \eqref{NM} with $L=N=M$ and $i = 0$. We get $\Hom_A(M,M) \simeq A$.
   It follows, again by \eqref{NM}, that
  \begin{equation} \label{NM2} \Ext^{i}_A(M,M)\simeq
    H_{*}^{i}(U,{\cH}om(\tilde{M},\tilde{M})) \simeq 
H^{i+1}_{JA}(\Hom_A(M,M)) = 0
\end{equation}
for $0 < i < j$, whence the result.
\end{proof}

\begin{remark} \label{SchSv} We consider the vanishing of $\Ext_A^i(M,M)$ in
  Theorem~\ref{Amodulethm} as mainly known (Schlessinger, see \cite{KleiLan},
  Prop.\! 2.2.3 for $i=1$) since it, as in \cite{S2} and \cite{Sv2}, is rather
  clear how to generalize \cite{KleiLan}, Prop.\! 2.2.3 to a non-smooth $X$
  and to $i > 1$ (e.g. \cite{Sv2}, Rem.\! 2.5). Our proof is, however, very
  short and uses more directly Grothendieck's long exact sequence of
  $\Ext$-groups appearing in \cite{SGA2}, exp.\! VI, from which \eqref{NM} is
  deduced.
 \end{remark}
 \begin{remark} \label{rile} It is clear from the proof that the theorem also
   holds for $c=1$. In this case we can only use the result for $j \le 2$
   because the largest possible value of $\depth_{I_{t-1}(\cA)A}A$ is $3$.
   Thus our proof implies the known rigidity of $M$ 
(\cite{KleiLan}, Prop.\! 2.2.3 and \cite{Ile}, Thm.\! 2). We
   continue to restrict ourselves to $c \ge 2$ and refer to \cite{Ile} for a
   nice study when $c=1$.
\end{remark}
In running some Macaulay 2 computations (\cite{Mac}) in the situation of
Theorem~\ref{Amodulethm} we were surprised to see that also $ \Ext_A^i(M,M)=0$
for $i=j$ in the examples. This observation led us to try to prove
Theorem~\ref{Amodulethm} under the assumption $\depth_{I_{t-1}(\cA)A}A \ge j$.
The natural case where this happens and where we succeed is as follows. Let
$B=R/I_t(\cB)$ and suppose $\depth_{JB}B \ge j+1$ with $J=I_{t-1}(\cB)$ where
${\cB}$ is obtained by deleting a column of ${\cA}$. Then since $I_{t}(\cA)
\subset I_{t-1}(\cB) \subset I_{t-1}(\cA)$, it follows that $\depth_{JA}A \ge
j$. Thus we may take $j= \codim_{Y}Sing(Y)-1= \dim X - \dim Sing(Y)$ in the
following theorem.

\begin{theorem} \label{Amodulethm2} Let $ B \to A$ be quotients of $R$ defined
  by the vanishing of the maximal minors of ${\cB}$ and ${\cA}$ respectively
  where ${\cB}$ is obtained by deleting some column of ${\cA}$. Let $M$ be the
  finitely generated graded module over $R$ of maximal grade defined by
  ${\cA}$, i.e. $M:=M_{\cA}$, cf. \eqref{Mi}. Let $j \ge 2$ be an integer and
  suppose $\depth_{I_{t-1}(\cB)B}B \ge j+1$. Then $\Hom_A(M,M) \simeq A$ and
 $$ \Ext_A^i(M,M) = 0 \  {\rm for} \ 1 \le i \le j-1 \ .$$
\end{theorem} 
\begin{proof} Since $M$ has maximal grade we get that $A$ and hence $B$
   are standard determinantal rings (by \cite{B} since we may suppose the
   matrix $\cA$ is minimal). It follows that $N:=M_{\cB}$
  has maximal grade and we can apply Theorem~\ref{Amodulethm} to $N$. We get $
  \Ext_B^i(N,N) = 0$ for $ 1 \le i \le j-1.$ Moreover we have $\Hom_A(M,M)
  \simeq A$ by the proof of Theorem~\ref{Amodulethm}.

  Now we consider the exact sequence
\begin{equation}\label{AMi}
0\longrightarrow B(-a_{t+c-2}) \longrightarrow N
\longrightarrow M \longrightarrow 0
\end{equation}
induced by \eqref{Mi} and we put $B_a:= B(-a_{t+c-2})$. We {\it claim} that $
\Ext_B^{1}(M,M) $ is isomorphic to $ \Hom_B(B_a,M) \simeq M(a_{t+c-2})$ and that $ \Ext_B^i(M,M) =
0$ for $ 2 \le i \le j-1.$ To see it we apply $\Hom_B(-,M)$ and $\Hom_B(N,-)$
to  \eqref{AMi}. Their long exact sequences
fit into the following diagram {\footnotesize
\begin{equation}\label{diagram}
\begin{array}{cccccccccccc}
   && && \downarrow  &&    \\ 
  & &  &  &  \Ext_B^i(N,N) &  &  & & &  \\ && &&
  \downarrow  &&     \\ \to \ \Ext_B^{i-1}(B_a,M) & \rightarrow &
  \Ext_B^i(M,M) &  \to & 
  \Ext_B^i(N,M) & \rightarrow & \Ext_B^{i}(B_a,M) & 
  \rightarrow & \Ext_B^{i+1}(M,M) \rightarrow
  \\ && & &  \downarrow  && \\
  & & & & \Ext_B^{i+1}(N,B_a) \\
  && & & \downarrow 
 \end{array}
\end{equation}}
where $ \Ext_B^{i}(B_a,M)=0$ for $i > 0$. To see that also
$\Ext_B^{i+1}(N,B_a)=0$ for $1 \le i+1 \le j$, we first notice that $
\Ext_B^{i+1}(N,B_a) \simeq  \Ext_B^{i+1}(N \otimes K_B,K_B(-a_{t+c-2}))$
(\cite{her}, Satz 1.2).
Moreover since $K_{B}(n+1) \simeq S_{c-2}N(\ell_{c-1})$, cf.\! \eqref{ell}, we
see that $N \otimes K_B$ and $ S_{c-1}N $ are closely related. Indeed up to
twist they are isomorphic if we restrict to $U_B:=\Proj(B)-V(I_{t-1}(\cB)B)$.
Hence if we let $\Lambda$ be the kernel of the natural surjective map $
S_{c-2}N \otimes_B N \to S_{c-1}N$, it follows that $\Supp{\Lambda} \subset
V(I_{t-1}(\cB)B)$, e.g. we get $ \Ext_B^{i+1}(\Lambda,K_B)=0$ for $i+1 \le j$
by the assumption $\depth_{I_{t-1}(\cB)B}B \ge j+1$. Now we recall that $
S_{c-1}N$ is a maximal CM $B$-module (\cite{eise}). It follows that $
\Ext_B^{i+1}( S_{c-1}N,K_B)=0$ for $i+1 >0$. Since the sequence
$$ \to  \Ext_B^{i+1}( S_{c-1}N,K_B) \to \Ext_B^{i+1}(  S_{c-2}N
\otimes_B N ,K_B) \to \Ext_B^{i+1}( \Lambda ,K_B) \to $$ is exact, we deduce
that $\Ext_B^{i+1}(N,B_a)=0$ for $0 \le i \le j-1$. Then using the big diagram
above we get the claim provided we can prove that the surjective map $
\Hom_B(B_a,M) \to \Ext_B^{1}(M,M) $ is an isomorphism. To prove it we continue
the horizontal sequence in the big diagram on the left hand side and we get
the exact sequence $$ 0 \rightarrow \Hom_B(M,M) \to \Hom_B(N,M) \rightarrow
\Hom_B(B_a,M) \rightarrow \Ext_B^{1}(M,M) \rightarrow 0 \ .$$ The two leftmost
$\Hom$-modules are easily seen to be isomorphic to $A$, e.g. $ \Hom_B(N,M)
\simeq 
H_{*}^{0}(U_B,{\cH}om(\tilde{N},\tilde{M})) \simeq A$ because $\tilde{N}
\otimes O_X \simeq \tilde{M}$ is invertible over $U_B \cap X$ and the {\it
  claim} is proved.

It remains to compare the groups $ \Ext_A^i(M,M)$ and $ \Ext_B^i(M,M)$ for
which we have a well-known spectral sequence $E_2^{p,q}:=
\Ext_A^p(\Tor_q^B(A,M),M)$, converging to $ \Ext_B^{p+q}(M,M)$, at our
disposal. Since $E_2^{p,0} \simeq \Ext_A^p(M,M)$ we must show 
\begin{equation} \label{getthm}
 E_2^{p,0} = 0
\ \ {\rm for} \ 1 \le p <j \ .
\end{equation}
Noticing that $I_{X/Y} \simeq N(a_{t+c-2})^*$ by \eqref{Di} in which $\tilde
N(a_{t+c-2})\arrowvert_{U_B}$ is an invertible $\cO_Y$-Module over $U_B
\subset Y$, we get that the sheafification of $ \Tor_q^B(A,M) \simeq
\Tor_{q-1}^B(I_{X/Y},M)$ vanishes over $U_B \cap X$ for $q \ge 2$. Hence
$\Supp { \Tor_q^B(A,M)} \subset V(I_{t-1}(\cB)A)$, and we get $E_2^{p,q}=0$
for $q \ge 2$ and $p < j$ by the assumption $\depth_{I_{t-1}(\cB)B}B \ge j+1$
which leads to $\depth_{I_{t-1}(\cB)A}A \ge j$. Note also that $E_2^{p,1} = 0$
for $0 < p < j-1$ because by \eqref{NM},
\begin{equation*} 
E_2^{p,1}=
\Ext_A^p(I_{X/Y} \otimes_B M,M) \simeq H_{*}^{p}(U,{\cH}om_{{\cO}_X}(\cI_{X/Y}
\otimes \tilde{M},\tilde{M})) \simeq H^{p+1}_{I_{t-1}(\cB)A}M(a_{t+c-2}) =0
\end{equation*}
where $U:=U_B \cap X$. Indeed $\cI_{X/Y}\otimes \tilde{M} \simeq \tilde
{N}(a_{t+c-2})^* \otimes_{\cO_Y} {\cO_X }\otimes \tilde{M} \simeq
\cO_X(-a_{t+c-2})$ over $U$. In the same way $$E_2^{0,1}= \Hom_A(I_{X/Y}
\otimes_B M,M) \simeq H_{*}^{0}(U,\tilde{M}(a_{t+c-2})) \simeq M(a_{t+c-2}) \
.$$ The spectral sequence leads therefore to an exact sequence 
\begin{equation} \label{specmid}
 0 \to
E_2^{1,0} \to \Ext_B^1(M,M) \to E_2^{0,1} \to E_2^{2,0} \to \Ext_B^2(M,M)\to
0
\end{equation}
and to isomorphisms $ E_2^{p,0} \simeq \Ext_B^p(M,M)$ for $2 < p <j$. We
already know $ \Ext_B^p(M,M)=0$ for $2 \le p <j$ by the proven claim. Hence we
get \eqref{getthm} provided we can show that the ``pushforward morphism''
$\Ext_B^1(M,M) \to E_2^{0,1} \simeq M(a_{t+c-2})$ of \eqref{specmid} is an
isomorphism. Since it is clear that this morphism is compatible with the
isomorphism $\Ext_B^1(M,M) \simeq \Hom_B(B_a,M) \simeq M(a_{t+c-2})$ which we
proven in the {\it claim} (using \eqref{AMi}), 
we get the theorem.
\end{proof}

\begin{remark} The proof of Theorem~\ref{Amodulethm2} even shows $
  \Ext_B^i(M,M) = 0 \ {\rm for} \ 2 \le i \le j-1 \ $.
 \end{remark}

 As we see, the proof of Theorem~\ref{Amodulethm2} is technically much more
 complicated that the proof of Theorem~\ref{Amodulethm}, even though we are
 only able to weaken the depth assumption on $A$ in some cases (namely when
 the two algebras in $R/I_{t-1}(\cB) \twoheadrightarrow R/I_{t-1}(\cA)$ have
 the same dimension). This improvement is, however, important in low
 dimensional cases in which the radical of ${I_{t-1}(\cB)}$ often satisfies
\begin{equation} \label{delsmooth}
\goth m = \sqrt{I_{t-1}(\cB)}
\end{equation}
and hence $ \sqrt{ I_{t-1}(\cB)} = \sqrt{I_{t-1}(\cA)}$. For short we say that
we get an l.c.i. scheme by deleting some column if \eqref{delsmooth} holds. We
immediately get from the theorems

\begin{corollary} \label{Amodulecor1} Let $X=\Proj(A)$, $A=R/I_{t}(\cA)$ be a
  standard determinantal scheme, let $M=M_{\cA}$ and suppose either
  $\depth_{I_{t-1}(\cA)A}A \ge 3$, or just $\dim X \ge 1$ provided we get an
  l.c.i. (e.g. a smooth) scheme by deleting some column of $\cA$. Then
  $\Hom_A(M,M) \simeq A$ and
 $$ \Ext_A^1(M,M) = 0 \ .$$ 
\end{corollary}

\begin{corollary} \label{Amodulecor2} Let $X=\Proj(A)$, $A=R/I_{t}(\cA)$ be a
  standard determinantal scheme, let $M=M_{\cA}$ and suppose $\dim X \ge 1$.
  Moreover suppose the polynomials $f_{ij}$ of degrees $a_j-b_i$ in a
  presentation matrix $(f_{ij})$ of $M$ are chosen general enough and suppose
  $a_{i-2}\ge b_{i}$ for $2 \le i \le t$. Then $\Hom_A(M,M)
  \simeq A$ and $ \ \Ext_A^1(M,M) = 0 \ .$
\end{corollary}
\begin{proof} We may suppose that the codimension of $X$ in $\proj{n}$ is $c
  \ge 3$ since $M$ is a twist of the canonical module of $A$ if $c=2$ in which
  case the conclusion is well known. Suppose $\dim X = 1$. Then
  Remark~\ref{dep} with $\alpha = 2$ shows that both $X=X_c$ and $Y:=X_{c-1}$
  are smooth because $X$ is general. If, however, $\dim X \ge 2$, then
  Remark~\ref{dep} still applies to $X=X_c$ and we get
  $\depth_{I_{t-1}(\cA)A}A \ge 3$. Hence in any case we conclude by
  Corollary~\ref{Amodulecor1}.
 \end{proof}

 In deformation theory it is important to know when $ \Ext_A^2(M,M)$
 vanishes.

 \begin{corollary} \label{Amodulecor3} Let $X=\Proj(A)$, $A=R/I_{t}(\cA)$ be a
   standard determinantal scheme, let $M=M_{\cA}$ and suppose either
   $\depth_{I_{t-1}(\cA)A}A \ge 4$, or just $\dim X \ge 2$ provided we get an
   l.c.i. (e.g. a smooth) scheme by deleting some column of $\cA$. Then
   $\Hom_A(M,M) \simeq A$ and
 $$ \Ext_A^i(M,M) = 0 \ \ {\rm \ for \ } i=1 \  {\rm and} \ 2 \ .$$ 
\end{corollary}
\begin{proof} This follows immediately from Theorem~\ref{Amodulethm2} and
  Theorem~\ref{Amodulethm}
\end{proof}
\begin{corollary} \label{Amodulecor4} Let $X=\Proj(A)$, $A=R/I_{t}(\cA)$ be a
  standard determinantal scheme, let $M=M_{\cA}$ and suppose $\dim X \ge 2$.
  Moreover suppose the polynomials $f_{ij}$ of degrees $a_j-b_i$ in a
  presentation matrix $(f_{ij})$ of $M$ are chosen general enough and suppose
  $a_{i-\min (3,t)}\ge b_{i}$ for $\min (3,t)\le i \le t$. Then
  $\Hom_A(M,M) \simeq A$ and
  $$ \Ext_A^i(M,M) = 0 \ \ {\rm \ for \ } i=1 \  {\rm and} \ 2 \ .$$ 
\end{corollary}
\begin{proof} We may again suppose that $c \ge 3$. Now if $\dim X = 2$, then
  Remark~\ref{dep} with $\alpha = 3$ shows that both $X=X_c$ and $Y:=X_{c-1}$
  are smooth. If, however, $\dim X \ge 3$, then Remark~\ref{dep} still applies
  to $X=X_c$ and we get $\depth_{I_{t-1}(\cA)A}A \ge 4$. Thus in any case we
  conclude by Corollary~\ref{Amodulecor3}. 
\end{proof}

\begin{remark} \label{Amodulerem} Also the results of this section admit {\it
    substantial generalizations} since the assumption that $R$ is a
  polynomial ring is not necessary. For instance let $R$ be any graded
  quotient of a polynomial ring $k[x_0, \dots ,x_n]$ with the standard grading
  where $k$ is any field.
  In Theorem~\ref{Amodulethm} it suffices to have $\depth_{I_{t-1}(\cA)A}M =
  \depth_{I_{t-1}(\cA)A}A $ and the depth assumption of that theorem to see
  that the proof works ($\tilde M \arrowvert_U$ locally free of rank one holds
  in general by \cite{BH}, Lem.\! 1.4.8). Moreover in
  Theorem~\ref{Amodulethm2}, Corollary~\ref{Amodulecor1} and
  Corollary~\ref{Amodulecor3} we use a few places that $R$ is Cohen-Macaulay
  in which case we get $\depth_{I_{t-1}(\cA)A}M = \depth_{I_{t-1}(\cA)A}A $ by
  \cite{eise}, Cor.\! A2.13. So all the mentioned results hold if $\Proj(R)$
  is any ACM-scheme (i.e. $R$ is CM). The remaining corollaries hold as well
  if $\Proj(R)$ is a smooth ACM scheme and $k = \overline k$ by
  Remark~\ref{dep}. Indeed Remark~\ref{dep} is really a result for
  determinantal subschemes of any smooth variety $W$, not only when $W=\PP^n$.
\end{remark} 

\section{deformations of modules and determinantal schemes}

The main goal of this section is to show a close relationship between the
local deformation functor, ${\rm Def}_{M/R}$, of the graded $R$-module $M =
M_{\cA}$ and the corresponding local functor, ${\rm Def}_{A/R}$, of deforming
the determinantal ring $A=R/{\rm ann}(M)$ as a graded quotient of $R$. We will
see that these functors are isomorphic (resp. the first is a natural
subfunctor of the other) provided $\dim X \ge 2$ (resp. $\dim X = 1$) and
$X=\Proj(A)$ is general. If $\dim X = 1$, the mentioned subfunctor is indeed
the functor which corresponds to deforming the determinantal $k$-algebra $A$
as a {\it determinantal} quotient of $R$ (Definition~\ref{subf}). Combining
with results of previous sections and the fact that ${\rm Def}_{A/R}$ is the
same as the local Hilbert (scheme) functor of $X$ if $\dim X \ge 1$ by
\eqref{Grad}, we get the main results of this paper; the dimension formula for
$ W(\underline{b};\underline{a})$ and the generically smoothness of
$\Hi^p(\proj{n})$ along $ W(\underline{b};\underline{a})$. The comparison is
mostly to understand well a spectral sequence comparing the tangent and
obstruction spaces of the mentioned deformation functors and to use the
theorems of the previous sections. This spectral sequence is also important in
R. Ile's PhD thesis \cite{I01}, and in his papers \cite{Ile} and \cite{IleTr}
(see Remark~\ref{ruile}). In the following we suppose
$A$ is generically a complete intersection ($\depth_{I_{t-1}(\cA)A}A \ge
1$), i.e. that $X=\Proj(A)$ is a good determinantal scheme.

Consider the well-known spectral sequence
$$E_2^{p,q}:=
\Ext_A^p(\Tor_q^R(A,M),M) \ \ \ \Rightarrow \ \ \ \Ext_R^{p+q}(M,M) \ ,$$ and
note that $E_2^{p,0} \simeq \Ext_A^p(M,M)$ and $ \Tor_q^R(A,M) \simeq
\Tor_{q-1}^R(I_{X},M)$ for $q \ge 1$. The spectral sequence leads to the
following exact sequence 
{\small
\begin{equation} \label{specseq} 0 \to \Ext_A^1(M,M) \to \Ext_R^1(M,M) \to
E_2^{0,1} \to \Ext_A^2(M,M) \to \Ext_R^2(M,M)\to E_2^{1,1} \to \ . 
\end{equation}
} Indeed $ E_2^{0,2}= \Hom(\Tor_{2}^R(A,M),M) = 0$ because $\Tor_{2}^R(A,M)$
is supported in $V({I_{t-1}(\cA)A})$. Moreover $$E_2^{0,1} \simeq \Hom_A(I_X
\otimes_R M,M) \simeq \Hom_R(I_X,\Hom_R(M,M)) \ ,$$ and see \cite{IleTr},
Def.\! 3 for an explicit description of $ \Ext_R^1(M,M) \to E_2^{0,1}$. In
our situation we
recall that $\depth_{I_{t-1}(\cA)A}A \ge 1$ lead to $\Hom_A(M,M) \simeq A$ by
\cite{KM}, Lem.\! 3.2. It follows that the edge homomorphism $ \Ext_R^1(M,M)
\to E_2^{0,1}$ of the spectral sequence above induces a natural map
\begin{equation} \label{specseq2} 
   _0\! \Ext_R^1(M,M) \longrightarrow \  (E_2^{0,2})_0 \simeq \ _0\!
   \Hom_R(I_X,A) 
\end{equation}
 between the tangent spaces of the two deformation functors ${\rm Def}_{M/R}$
 and ${\rm Def}_{A/R}$ respectively.
 Even though we only partially use the spectral sequence in the proof below,
 Theorem~\ref{compthm} is fully motivated by the spectral sequence.
 \begin{definition} \label{subf} Let $X=\Proj(A)$, $A=R/I_t(\cA)$, be a good
   determinantal scheme and let $\underline {\ell}$ be the category of
   artinian $k$-algebras (cf. the text before \eqref{al}). Then the local
   deformation functor ${\rm Def}_{A \in W(\underline{b};\underline{a})}$,
   defined on $\underline {\ell}$, is the subfunctor of ${\rm Def}_{A/R}$
   given by: \vspace{0.2cm}
\begin{equation*}
{\rm Def}_{A \in
  W(\underline{b};\underline{a})}(T) = \left\{ A_T \in {\rm
    Def}_{A/R}(T) \arrowvert A_T=R_T/I_t(\cA_T) {\rm \ for \ some \
    matrix \     \cA_T \ lifting \ \cA \ to \ } T \right\}.
\end{equation*}
\vspace{-0.3cm}
\end{definition}

Note that there is a natural map \ ${\rm Def}_{M/R} \to {\rm Def}_{A \in
  W(\underline{b};\underline{a})}$ because for every graded deformation $M_T$
of $M$ to $T$ there exists a matrix $\cA_{T}$ whose induced morphism has
$M_{{T}}$ as cokernel (see the first part of the proof of
Theorem~\ref{modulethm}) and because different matrices inducing the same
$M_T$ define the same ideal of maximal minors (Fittings lemma, \cite{eise},
Cor.\! 20.4). The map is surjective since we can use the matrix $\cA_T$ in
Definition~\ref{subf} to define $M_{{T}} \in {\rm Def}_{M/R}(T)$.

The phrase ``for some matrix $ \cA_T$ lifting $\cA$ to $T$'' which means that
there exists a homogeneous matrix $ \cA_T$ lifting $\cA$ to $T$, may be
insufficient for forcing ${\rm Def}_{A \in W(\underline{b};\underline{a})}$ to
have nice properties. For instance we do not know whether ${\rm Def}_{A \in
  W(\underline{b};\underline{a})}$ is pro-representable, or even has a hull,
since we have no proof for the surjectivity of
\begin{equation} \label{schl}
{\rm Def}_{A \in
    W(\underline{b};\underline{a})}(T_1 \times_S T_2) \longrightarrow {\rm
    Def}_{A \in W(\underline{b};\underline{a})}(T_1) \times_{ {\rm Def}_{A \in
      W(\underline{b};\underline{a})}(S)} {\rm Def}_{A \in
    W(\underline{b};\underline{a})}(T_2) \ 
  \end{equation}
  for every pair of morphisms $T_i \to S$, $i=1,2$, in $\underline {\ell}$
  with $T_2 \twoheadrightarrow S$ small (see Schlessinger's main theorem in
  \cite{S}). In \cite{Sha} Schaps solves a related problem by
  assuming that $\cA$ has the unique lifting property and she gets some results
  on the existence of a hull for determinantal non-embedded deformations.
  In our context, assuming $ _0\! \Ext_A^1(M,M)= 0$, then we shall
  see that \ ${\rm Def}_{A \in W(\underline{b};\underline{a})}$ behaves well
  because for every element of \ ${\rm Def}_{A \in
    W(\underline{b};\underline{a})}(T)$ there exists a {\em unique} module
  $M_{\cA_T}$ even though $\cA_T$ is not unique. 

  Indeed let $D:=k[\epsilon]/(\epsilon^2)$ be the dual numbers and
  let $$\lambda := \dim \ _0\! \Ext_R^1(M,M) = \lambda_c + K_3+K_4+...+K_c \
  ,$$ cf.\! Theorem~\ref{modulethm}. Recalling that $
  W(\underline{b};\underline{a})$ is a certain quotient of an open irreducible
  set in the affine scheme $\VV =\Hom_{{\mathcal O}_{\PP^{n}}}({\mathcal
    G^*},{\mathcal F^*})$
  whose rational points correspond to $t \times (t+c-1)$ matrices and that
  $\dim W(\underline{b};\underline{a}) \le \lambda$ (\cite{KM}, p.\! 2877 and
  Thm.\! 3.5), we
  get

  \begin{theorem} \label{compthm} Let $X=\Proj(A)$, $A=R/I_{t}(\cA)$ be a good
    determinantal scheme. If \ $ _0\! \Ext_A^1(M,M)= 0$ then the functor ${\rm
      Def}_{A \in W(\underline{b};\underline{a})}$ is pro-representable, the
    pro-representing object has dimension $\dim
    W(\underline{b};\underline{a})$ and $${\rm Def}_{A \in
      W(\underline{b};\underline{a})} \simeq {\rm Def}_{M/R} \ . $$ Hence
    ${\rm Def}_{A \in W(\underline{b};\underline{a})}$ is formally smooth.
    Moreover the tangent space of \ ${\rm Def}_{A \in
      W(\underline{b};\underline{a})}$ is the subvector space of $ _0\!
    \Hom_R(I_X,A)$ which corresponds to graded deformation $R_D \to A_D$ of $R
    \to A$ to $D$ of the form $A_D=R_D/I_t(\cA_D)$ for some matrix $\cA_D$
    which lifts $\cA$ to $D$. \vspace{0.15cm}

    If in addition \ $ _0\! \Ext_A^2(M,M)= 0$, then $ {\rm Def}_{M/R} \simeq
    {\rm Def}_{A \in W(\underline{b};\underline{a})} \simeq {\rm Def}_{A/R}$ and
    \ ${\rm Def}_{A/R}$ is formally smooth. Moreover every deformation of
    $X$ 
    comes from deforming $\cA$ (cf.\! Definition~\ref{everydef}).
  \end{theorem}

  \begin{proof} We already know that \ ${\rm Def}_{M/R}(T) \to {\rm Def}_{A
      \in W(\underline{b};\underline{a})}(T)$ is well defined and surjective.
    To see that it is injective, we will construct an inverse. Suppose
    therefore that there are two matrices $ (\cA_T)_1$ and $ (\cA_T)_2$
    lifting $\cA$ to $T$ and such that $I_t((\cA_T)_1) = I_T((\cA_T)_2)$. The
    two matrices define two graded deformations $M_1$ and $M_2$ of the
    $R$-modules $M$ to $R_T$ by Lemma ~\ref{defalpha}. Since, however, the two
    matrices define the {\em same} graded deformation $A_T:=R_T/
    I_T((\cA_T)_1)$ of $A$ to $T$, we get that $M_1$ and $M_2$ are two graded
    deformations of the $A$-module $M$ to $A_T$! Due to \ $ _0\!
    \Ext_A^1(M,M)= 0$, $\Hom_A(M,M) \simeq A$ and deformation theory, we
    conclude that $M_1 = M_2$ up to multiplication with a unit of $T$, i.e. we
    get a well defined map which clearly is an inverse.

    Since we have ${\rm Def}_{A \in W(\underline{b};\underline{a})} \simeq
    {\rm Def}_{M/R}$ and we know that $ {\rm Def}_{M/R}$ has a hull
    (\cite{siq}), it follows that ${\rm Def}_{A \in
      W(\underline{b};\underline{a})} $ has a hull (or one may easily show the
    surjectivity of \eqref{schl} directly by using the uniqueness of
    $M_{\cA_T}$). Note that the injectivity of \eqref{schl} follows from ${\rm
      Def}_{A \in W(\underline{b};\underline{a})}$ being a subfunctor of the
    pro-representable functor ${\rm Def}_{A/R}$ (\cite{K04}, Prop.\! 9),
    whence ${\rm Def}_{A \in W(\underline{b};\underline{a})}$ is
    pro-representable by \cite{S}. Moreover using ${\rm Def}_{A \in
      W(\underline{b};\underline{a})} \simeq {\rm Def}_{M/R}$ and
    Theorem~\ref{modulethm} we get that ${\rm Def}_{A \in
      W(\underline{b};\underline{a})} $ is formally smooth and that $\dim H =
    \lambda$ where $H$ is the pro-representing object of ${\rm Def}_{A \in
      W(\underline{b};\underline{a})}$. The description of its tangent space
    follows from Definition~\ref{subf} and \eqref{specseq}\!
    -\eqref{specseq2} since $ _0\! \Hom_R(I_X,A)$ is the tangent space of $
    {\rm Def}_{A/R}$.

    So far we know $\dim W(\underline{b};\underline{a}) \le \lambda = \dim H$.
    To see that $\dim H = \dim W(\underline{b};\underline{a})$, it
    suffices to see that the family of determinantal rings over $H$,
    corresponding to the ``universal object'' of ${\rm Def}_{A \in
      W(\underline{b};\underline{a})}$, is
    algebraizable. 
    This is clear in our context, (see the explicit description of $H$ in the
    proof of \cite{L}, Thm.\! 4.2.4). 
    Indeed take $\lambda$ independent elements of $\ _0\! \Ext_R^1(M,M)\simeq
    {\rm Def}_{A \in W(\underline{b};\underline{a})}(D)$, let $\cA+ \epsilon
    \cA_1,...,\cA+ \epsilon \cA_{\lambda}$ be corresponding presentation
    matrices of the elements (i.e. modules), and let $\cA_T := \cA+ t_1
    \cA_1+...+t_{\lambda}\cA_{\lambda}$ (linear combination in the parameters
    $t_{k}$) where $T$ be the polynomial ring $T = k[t_1,...,t_{\lambda}]$.
    Then the algebraic family $A_T:=R_T/I_t(\cA_T)$ is $T$-flat at $(0,...,0)
    \in \Spec(T)$ (Lemma~\ref{defalpha}) and hence flat in a neighborhood and
    we get what we want.

    Finally we suppose $ _0\! \Ext_A^2(M,M)= 0$. Using ${\rm Def}_{A \in
      W(\underline{b};\underline{a})}(D) \simeq \ _0\! \Ext_R^1(M,M)$ and the
    spectral sequence \eqref{specseq} we get isomorphisms ${\rm Def}_{A \in
      W(\underline{b};\underline{a})}(D) \simeq \ _0\! \Hom_R(I_X,A) \simeq
    {\rm Def}_{A/R}(D)$ of tangent spaces. To show ${\rm Def}_{A \in
      W(\underline{b};\underline{a})}(T) \simeq {\rm Def}_{A/R}(T)$ for any
    $(T,\mathfrak{m}_T)$ in $\underline{\ell}$, we may by induction suppose
    $\mathfrak{m}_T^{r+1}=0$ and ${\rm Def}_{A \in
      W(\underline{b};\underline{a})}(T/\mathfrak{m}_T^{r}) \simeq {\rm
      Def}_{A/R}(T/\mathfrak{m}_T^{r})$. Consider the commutative diagram
\begin{equation}
\begin{array}{cccccccccccc}
 &   &  & 
  {\rm Def}_{A \in
    W(\underline{b};\underline{a})}(T) & \hookrightarrow &  {\rm Def}_{A/R}(T) & 
  \\ && & \downarrow  &&  \downarrow  && \\
  & & & {\rm Def}_{A \in
    W(\underline{b};\underline{a})}(T/\mathfrak{m}_T^{r}) &\simeq &  {\rm
    Def}_{A/R}(T/\mathfrak{m}_T^{r})
 \end{array}
\end{equation}
and notice that the leftmost vertical map is surjective since $ {\rm Def}_{A
  \in W(\underline{b};\underline{a})}$ is formally smooth. Hence for a given
$A_T \in {\rm Def}_{A/R}(T)$ there exists $A'_T \in {\rm Def}_{A/R}(T)$ such
that $A'_T \simeq R_T/I_t(\cA_T)$ for some matrix $\cA_T$ which lifts a matrix
$\cA_{T/\mathfrak{m}_T^{r}}$ defining $A_T \otimes_T T/\mathfrak{m}_T^{r}$ to
$T$. The difference of $A_T$ and $A'_T$ belongs to $ {\rm
  Def}_{A/R}(D)\otimes_k \mathfrak{m}_T^{r} \simeq {\rm Def}_{A \in
  W(\underline{b};\underline{a})}(D) \otimes_k \mathfrak{m}_T^{r}$, and
``adding'' it to $A'_T$ we get that $A_T \in {\rm Def}_{A \in
  W(\underline{b};\underline{a})}(T)$, whence ${\rm Def}_{A \in
  W(\underline{b};\underline{a})}(T) \simeq {\rm Def}_{A/R}(T)$. It follows
that the completion of the local ring $\cO_{\Hi ,(X)}$ of $ \Hi ^p(\PP^{n})$ at
$(X)$ is isomorphic to $H$. Since we in the preceding paragraph explicitly
constructed an algebraic determinantal family over some neighborhood of
$(0,...,0)$ in $ \Spec(k[t_1,...,t_{\lambda}])$ (thinking about it, we must
have $\cO_{\Hi,(X)} \simeq k[t_1,...,t_{\lambda}]_{(t_1,...,t_{\lambda})}$
since $k= \overline k$), we get that ``every deformation of $X$ comes from
deforming $\cA$`` and we are done.
\end{proof}
\begin{remark} \label{comprem} Let us endow the closed subset $ \overline{
    W(\underline{b};\underline{a})}$ of $ \Hi ^p(\PP^{n})$ with the reduced
  scheme structure (this is natural since ``the part $
  W(\underline{b};\underline{a})$ of $ \Hi ^p(\PP^{n})$ is unobstructed'' by
  the proof of Lemma~\ref{unobst}). Let $X=\Proj(A)$, $A=R/I_{t}(\cA)$ belong
  to $ W(\underline{b};\underline{a})$. Then the proof related to $\dim
  W(\underline{b};\underline{a}) = \lambda$ above imply that the Zariski
  tangent space, $( \mathfrak{m}_W/ \mathfrak{m}_W^2)^\vee$, of $
  W(\underline{b};\underline{a})$ satisfies
  \begin{equation} \label{dimWvee}
    ( \mathfrak{m}_W/ \mathfrak{m}_W^2)^\vee \ =  \ {\rm Def}_{A \in
      W(\underline{b};\underline{a})}(D) \ .
  \end{equation}
  In the proof we used $\dim W(\underline{b};\underline{a}) \le \lambda$
  (\cite{KM}, Thm.\! 3.5) to show $\dim W(\underline{b};\underline{a}) =
  \lambda$. We will now explain this inequality by a direct argument. Indeed
  take any $(X') \in W(\underline{b};\underline{a})$. Then there is a matrix
  $t \times (t+c-1)$ matrix $\cA'$ whose maximal minors define $X'$. By
  Lemma~\ref{defalpha} the matrix $\cA+x(\cA'-\cA)$, $x$ a parameter, defines
  a flat family of good determinantal schemes over some open set $U \subset
  \Spec(k[x]) \simeq \AAA^{1}$ containing $x=0$ and $x=1$. Thus to any $(X')
  \in W(\underline{b};\underline{a})$ there is a tangent direction, i.e.\! an
  element $ t_X$ of $( \mathfrak{m}_W/ \mathfrak{m}_W^2)^\vee \subset \ _0\!
  \Hom_R(I_X,A) = {\rm Def}_{A/R}(D)$, given by the matrix $\cA'-\cA$. By
  Definition~\ref{subf}, $ t_X \in {\rm Def}_{A \in
    W(\underline{b};\underline{a})}(D)$, thus $ ( \mathfrak{m}_W/
  \mathfrak{m}_W^2)^\vee \subset {\rm Def}_{A \in
    W(\underline{b};\underline{a})}(D)$ by the relationship between $
  W(\underline{b};\underline{a})$ and its Zariski tangent space. Taking
  dimensions we have shown $\dim W(\underline{b};\underline{a}) \le
  \lambda$. 
  Then the proof of Theorem~\ref{compthm} implies $\dim
  W(\underline{b};\underline{a}) = \lambda$ and hence we get \eqref{dimWvee}.
 \end{remark}
 \begin{remark} \label{ruile} If the assumption \ $ _0\! \Ext_A^1(M,M)= 0$ of
   Theorem~\ref{compthm} is not satisfied, then the local deformation
 functor $ {\rm Def}_{M/A}$ of deforming $M$ {\em as a graded $A$-module} and
 its connection to $ {\rm Def}_{M/R}$ may be quite complicated, see
 \cite{IleTr} which compares the corresponding non-graded functors using
 \eqref{specseq}. However, by the results of the preceding section, $
 \Ext_A^1(M,M)= 0$ and $\Hom_A(M,M) \simeq A$ are week assumptions for modules
 of maximal grade.
 \end{remark}
 We now deduce the main theorems of the paper. In the first theorem we
 let $$\ext^2(M,M):= \dim \ker (\ _0\! \Ext_A^2(M,M) \to \ _0\!
 \Ext_R^2(M,M)\,)\, , \ \ {\rm cf. \ \eqref{specseq}} , $$ and {\em notice}
 that we write $ \Hi (\PP^{n})$ for $ \Hi ^p(\PP^{n})$ (resp. $\GradAlg(H)$)
 if $n-c \ge 1$ (resp. $n-c=0$), cf. the text accompanying \eqref{Grad} for
 explanations and notations.
 \begin{theorem} \label{Amodulethm3} Let $X=\Proj(A) \subset \PP^{n}$,
   $A=R/I_{t}(\cA)$ be a good determinantal scheme of $
   W(\underline{b};\underline{a})$ of dimension $n-c \ge 0$, let $M=M_{\cA}$
   and suppose \ $ _0\! \Ext_A^1(M,M) = 0$. Then
 $$\dim   W(\underline{b};\underline{a}) =  \lambda_c + K_3+K_4+...+K_c \ .
 $$ Moreover, for the codimension of $
 W(\underline{b};\underline{a})$ in $ \Hi (\PP^{n})$ in a neighborhood of
 $(X)$ we have
\begin{equation*} \label{codims} \dim_{(X)} \Hi (\PP^{n}) - \dim
  W(\underline{b};\underline{a}) \le \ext^2(M,M) \ ,
 \end{equation*}
 with equality if and only if $ \Hi (\PP^{n})$ is smooth at $(X)$. In
 particular these conclusions hold if $\depth_{I_{t-1}(\cA)A}A \ge 3$, {\rm
   \underline{or}} if $n-c \ge 1$ and we get an l.c.i. (e.g. a smooth)
 scheme by deleting some column of $\cA$.
\end{theorem}
\begin{proof} This follows from Theorem~\ref{compthm},
  Theorem~\ref{modulethm}, \eqref{specseq}\! -\eqref{specseq2} and
  Corollary~\ref{Amodulecor1}.
\end{proof}
\begin{corollary} \label{Amodulethm4} Given integers $a_0\le a_1\le ... \le
  a_{t+c-2}$ and $b_1\le ...\le b_t$, we suppose $n-c\ge 1$  and
   $a_{i-2}-b_i \ge 0$ for $2 \le i \le t$. Then
 $$\dim   W(\underline{b};\underline{a}) =  \lambda_c + K_3+K_4+...+K_c \ .$$
 provided $ \dim W(\underline{b};\underline{a}) \ne \emptyset$. In particular
 Conjecture 4.1 of \cite{KM09} holds in the case $ n-c \ge 1$.
\end{corollary}
\begin{proof} This follows from Theorem~\ref{Amodulethm3} and
  Corollary~\ref{Amodulecor2} since Conjecture 4.1 of \cite{KM09} is
  Conjecture~\ref{conjnew1} of this paper (and remember that we always suppose
  $c \ge 2$ and $t \ge 2$).
\end{proof}

\begin{remark} Even for zero-dimensional determinantal schemes ($n-c = 0$) the
  assumption \ $_0\! \Ext_A^1(M,M)= 0$ seems very week, and hence we almost
  always have the conjectured value of $\dim W(\underline{b};\underline{a})$.
  Thus Theorem~\ref{Amodulethm3} completes Theorem 4.19 of \cite{K09} in the
  zero-dimensional case. Indeed in computing many examples by Macaulay 2 we
  have so far only found \ $_0\! \Ext_A^1(M,M) \ne 0$ for examples outside the
  range of Conjecture~\ref{conjnew1}.
 \end{remark}
 Note that $ \overline{ W(\underline{b};\underline{a})}$ is not always an
 irreducible component of $ \Hi (\PP^{n})$. An example showing this was given
 in \cite{KMMNP}, Ex.\! 10.5, and many more were found in \cite{K09}, Ex.\!
 4.1, in which there are examples for every $c \ge 3$ (the matrix
 is linear except for the last column). All examples satisfy $n-c \le 1$.
 Indeed \cite{K09} contains exact formulas for the codimension of $ \overline{
   W(\underline{b};\underline{a})}$ in $ \Hi (\PP^{n})$ under some
 assumptions. Further investigations in \cite{KM09} led us to conjecture that
 $ \overline{ W(\underline{b};\underline{a})}$ is an irreducible component
 provided $n-c \ge 2$. Now we can prove it!

 \begin{theorem} \label{Amodulethm5} Let $X=\Proj(A) \subset \PP^n$,
   $A=R/I_{t}(\cA)$ be a good determinantal scheme of $
   W(\underline{b};\underline{a})$ of dimension $n-c \ge 1$, let $M=M_{\cA}$
   and suppose \ $ _0\! \Ext_A^i(M,M) = 0$ for $i=1$ and $2$. Then the Hilbert
   scheme $\Hi ^p(\PP^{n})$ is smooth at $(X)$,
 $$\dim_{(X)} \Hi ^p(\PP^{n}) =  \lambda_c + K_3+K_4+...+K_c \ , $$
 and every deformation of $X$ comes from deforming $\cA$. In particular this
 conclusion holds if $\depth_{I_{t-1}(\cA)A}A \ge 4$, {\rm \underline{or}} if
  $n-c \ge 2$ and we get an l.c.i. (e.g. a smooth) scheme by deleting some
 column of $\cA$.
\end{theorem}
\begin{proof} This follows immediately from Theorem~\ref{compthm},
  Theorem~\ref{modulethm} and Corollary~\ref{Amodulecor3}.
\end{proof}
\begin{corollary} \label{Amodulethm6} Given integers $a_0\le a_1\le ... \le
  a_{t+c-2}$ and $b_1\le ...\le b_t$, we suppose $n-c\ge 1$, $a_{i-2}-b_i \ge
  0$ for $2 \le i \le t$ and $ _0\! \Ext_A^2(M,M) = 0$ for a general
  $X=\Proj(A)$ of $ W(\underline{b};\underline{a})$. Then the closure
  $\overline{ W(\underline{b};\underline{a})}$ is a generically smooth
  irreducible component of the Hilbert scheme $\Hi ^p(\PP^{n})$ of
  dimension $$ \lambda_c + K_3+K_4+...+K_c \ . $$ In particular this
  conclusion holds if $n-c\ge 2$, $a_{i-\min (3,t)}\ge b_{i}$ for $\min
  (3,t)\le i \le t$ and $ \dim W(\underline{b};\underline{a}) \ne \emptyset$.
  It follows that Conjecture 4.2 of \cite{KM09} holds.
\end{corollary}
\begin{proof} This follows from Corollary~\ref{Amodulecor2},
  Theorem~\ref{Amodulethm5}, Lemma~\ref{unobst} and
  Corollary~\ref{Amodulecor4}, and note that Conjecture 4.2 of \cite{KM09} is
  the same as Conjecture~\ref{conjcomp} of this paper.
\end{proof}
Even in the one-dimensional case ($n-c = 1$) the assumption \ $_0\!
\Ext_A^2(M,M)= 0$ seems rather week, and we can often conclude as in
Corollary~\ref{Amodulethm6}. Note that if \ $_0\! \Ext_A^2(M,M)= 0$ for a
general $X$ of $ W(\underline{b};\underline{a})$ and $a_{i-2}-b_i \ge 0$ for $2
\le i \le t$, we get
\begin{equation} \label{homex}
 \ _0\!
\hom_R(I_X,A) = \ \lambda_c + K_3+K_4+...+K_c \ 
\end{equation}
by Corollary~\ref{Amodulecor2} and \eqref{specseq}\! -\eqref{specseq2}. So one
may alternatively compute $ \ _0\! \hom_R(I_X,A)$ and check if \eqref{homex}
holds, to conclude as in Corollary~\ref{Amodulethm6}
(cf. Theorem~\ref{compthm}). In \cite{K09}, Prop.\! 4.15 (which
generalizes \cite{KMMNP}, Cor.\! 10.15) we gave several criteria for
determining $\overline{ W(\underline{b};\underline{a})}$ in the
one-dimensional case. None of them apply in the following example.
\begin{example}  [determinantal curves in $\PP^{4}$, i.e. with $c =
  3$] \label{excu}  \ \label{va}

  Let $\cA = (f_{ij})$ be a $2 \times 4$ matrix whose entries are general
  polynomials of the same degree $f_{ij}=2$. The vanishing of all $2 \times 2$
  minors of $\cA$ defines a smooth curve $X$ of degree $32$ and genus $65$ in
  $\PP^{4}$. A Macaulay 2 computation shows \ $_0\! \Ext_A^2(M,M)= 0$. It
  follows from Corollary~\ref{Amodulethm6} that $\overline{
    W(\underline{b};\underline{a})}$ is a generically smooth irreducible
  component of $\Hi ^p(\PP^{4})$ of dimension $\lambda_3 + K_3 = 101$.

  Note that our previous method was to delete a column to get a matrix $\cB$
  and an algebra $B:=R/J$, $J:=I_{t}(\cB)$ and to verify $ _0\!
  \Ext^1_{B}(J/J^2,I/J)=0$ with $I=I_{t}(\cA)$. However, by Macaulay 2, $ _0\!
  \Ext^1_{B}(J/J^2,I/J)$ as well as $ _0\! \Ext^1_{A}(I/I^2,A)$, are
  $5$-dimensional and the approach of using Prop.\! 4.15 (i) does not apply
  (since $ _0\! \Ext^1_{B}(I/J,A) \ne 0$), neither do Prop.\! 4.15 (ii) nor
  (iii).
\end{example}

It is known that the vanishing of the cohomology group $H^1({\mathcal N}_X)$
(resp. $ \Ext^1_A(I_{X}/I_{X}^2,A)$) of a locally (resp. generically) complete
intersection $X \hookrightarrow \PP^{}$ implies that $X$ is unobstructed, and
that the converse is not true, e.g. we may have $H^1({\mathcal N}_X) \ne 0$
for $X$ unobstructed. Since we by Theorem~\ref{Amodulethm5} get that $X$ is
unobstructed by mainly assuming $n-c \ge 2$, one may wonder if we can prove a
little more, namely $H^1({\mathcal N}_X) = 0$. Indeed we can if $n-c \ge 3$.
More precisely recalling $\depth_{J}A = \dim A - \dim A/J$ we have

\begin{theorem} \label{normal} Let $X=\Proj(A) \subset \PP^n$,
  $A=R/I_{t}(\cA)$ be a standard determinantal scheme.
 
  {\rm (i)} If $\depth_{I_{t-1}(\cA)A}A \ge 4$ or equivalently, $\dim X \ge
  3 + \dim R/I_{t-1}(\cA)$, then $$ \Ext^i_A(I_{X}/I_{X}^2,A)=0 \ \ {\rm for}
  \ 1 \le i \le \dim X - 2 - \dim R/I_{t-1}(\cA), \ \ {\rm and
  }$$ $$H^i({\mathcal N}_X(v))=0 \ \ {\rm for} \ \ 1 \le i \le \dim X-2 \ \
  {\rm and \ every} \ v.$$

  {\rm (ii)} In particular if $\dim X \ge 3$,
  $a_{i-\min (3,t)}\ge b_{i}$ for $\min (3,t)\le i \le t$ and $X$ is general
  in $ W(\underline{b};\underline{a})$, then conclusions of  {\rm (i)}
  hold. If furthermore  $a_j \ge b_{i}$ for every $j$ and $i$, then
  $$ \Ext^i_A(I_{X}/I_{X}^2,A)=0 \ \ {\rm for}
  \ 1 \le i \le \min\{\dim X - 2, c-1\} \ .$$
\end{theorem}

\begin{proof} (i) Using \eqref{specseq}\! -\eqref{specseq2} and
  Corollary~\ref{Amodulecor3} 
  we get $ \Ext_R^1(M,M) \simeq \Hom_R(I_X,A)$. It follows that
  \begin{equation} \label{de} \depth \Hom_R(I_X,A) \ge \dim X 
  \end{equation} by Theorem~\ref{modulethm}. Thus
  the local cohomology group $H^{i}_{\goth m}(\Hom_R(I_X,A))$ vanishes for $i
  < \dim X$. Recalling that the sheafification of $ \Hom_R(I_X,A) \simeq
  \Hom_A(I_X/I_X^2,A)$ is 
  ${\mathcal N}_X$, we get $H_{*}^{i}({\mathcal N}_X)=0$ for $1 \le i < \dim X
  -1$,   whence we have the second vanishing of (i).

  Next let $r:=\depth_JA -1$ where $J$ is the ideal $I_{t-1}(\cA)A$ of $A$. It
  is known that \eqref{de} also implies $ \depth_{J} \Hom_R(I_X,A) \ge r$
  (e.g.\! \cite{KP2}, Lem.\! 28). Thus the local cohomology group
  $H^{i}_{J}(\Hom_A(I_X/I_{X}^2,A))$ vanishes for $i < r$ and we get the first
  vanishing of (i) by \eqref{NM} (letting $N=I_X/I_{X}^2$ and $L=A$).

  (ii) Finally we use Remark~\ref{dep} with $\alpha = 3$ (resp. $\alpha \ge
  (c+3)/2$) to see that $\depth_JA \ge 4$ (resp. $\codim_X V(J) \ge c+2$). In
  particular (i) applies to get the first statement. For the final statement,
  we recall the well known fact that $c+2$ is the largest possible value of
  the height of $J$ in $A$, whence $\codim_X V(J)= c+2$ with the usual
  interpretation that $c+2 = \codim_X V(J) > \dim X$ implies $V(J) =
  \emptyset$. This implies the theorem.
\end{proof}

For the algebra cohomology groups ${\rm H}^i(R,A,A)$ of Andr\'{e}-Quillen (cf.
\cite{AND}) 
we deduce

\begin{corollary} \label{Svanes} Let $A=R/I_{t}(\cA)$ be a standard
  determinantal graded $k$-algebra.
 
  {\rm (\! i)} If $\depth_{I_{t-1}(\cA)A}A \ge 4$ then
  \begin{equation*} \label{algcoh} {\rm H}^i(R,A,A) = 0 \ \ \ {\rm for} \ \ 2
    \le i \le \depth_{I_{t-1}(\cA)A}A -2\, .
  \end{equation*} 
  {\rm (\! ii)} If $\dim A \ge 4$, $a_j \ge b_{i}$ for every $j, i$ and
  $\Proj(A)$ is general in $ W(\underline{b};\underline{a})$, then
 \begin{equation*} {\rm H}^i(R,A,A) = 0 \ \ \ {\rm for}
 \ \ 2 \le i \le \min\{\dim A - 2, c\} \, .
  \end{equation*} 
\end{corollary}
\begin{proof} The spectral sequence relating algebra cohomology to
  algebra homology (\cite{AND}, Prop\! 16.1 or \cite{L}), implies, under the
  sole assumption $\depth_{I_{t-1}(\cA)A}A \ge 1$, that $$
  \Ext^i_A(I_{X}/I_{X}^2,A) \simeq {\rm H}^{i+1}(R,A,A) \, . $$ 
\end{proof}
\begin{remark} \label{normalrem} (i) The vanishing of $H_{*}^{i}({\mathcal
    N}_X)$ of Theorem~\ref{normal} is known if $c=3$ (\cite{KP2}, Lem.\! 35)
  or $c=4$ (\cite{KM}, Cor.\! 5.5). It $c=2$ even more is true by \cite{Bu}
  (or see \cite{KMMNP}, Cor.\! 6.5).

  (ii) Note that Corollary~\ref{Svanes} for so-called {\em generic}
  determinantal schemes is proved by Svanes (see \cite{b-v}, Thm.\! 15.10)
  while \cite{BC}, (1.4.3) shows the corollary for some non-generic
  determinantal schemes as well.

  (iii) As for $c=2$ one may hope that $H_{*}^{i}({\mathcal N}_X)=0$
  also for $i=\dim X - 1$. This is not true, as one may see through examples,
  using e.g. Macaulay 2. We have checked it for some surfaces in the range $3
  \le c < 6$ and always found it to be non-zero (cf. \cite{b-v}, 15.11).
\end{remark} 
\begin{remark} \label{Amodulerem2} In proving Theorem~\ref{normal} we used
  Corollary~\ref{Amodulecor4} to see that not only \ $ _0\! \Ext_A^i(M,M)$
  vanishes for $i=1, 2$, but in fact that the whole \ $ \Ext_A^i(M,M)$-group
  vanishes for $i=1$ and $2$. Arguing as in Theorem~\ref{compthm} and using
  the vanishing of the whole \ $ \Ext_A^i(M,M)$-group for $i=1$ and $2$, we
  may see that the non-graded deformation functor;
  \begin{equation*} {\rm Def}_{A/R}^{\rm non-gr}(T) = \left\{ R_T \to A_T
      \arrowvert A_T \ {\rm is \ } T{- \rm flat} \ {\rm and} \ \ A_T \otimes_T
      k \simeq A \right\}
\end{equation*}
in which a deformation $A_T$ of $R \to A$ to an artinian $T$ in $\underline
{\ell}$ is possibly non-graded,
is formally smooth provided $\depth_{I_{t-1}(\cA)A}A \ge 4$, {\rm
  \underline{or}} $\dim X \ge 2$ and we get an l.c.i. scheme by deleting some
column of $\cA$. This result is the best possible with regard to $\dim X \ge
2$ because one knows that ${\rm Def}_{A/R}^{\rm non-gr}$ is non-smooth for a
one-dimensional rational normal scroll $\Proj(A) \subset \PP^n$ for $n \ge 4$
(\cite{Pi}). Note also that we may deduce the result above for generic
determinantal schemes satisfying $\dim X \ge 3$ by works of Svanes
(\cite{b-v}, Thm.\! 15.10)
\end{remark} 

\begin{remark} \label{Amodulerem3} The results so far of this section admit
  {\it substantial generalizations} with respect to $R$ being a polynomial
  ring. Indeed we may let $R$ be any graded CM quotient of a polynomial ring
  $k[x_0, \dots ,x_n]$, $k = \overline k$, with the standard grading provided
  we in all results replace $\PP^{n}$ by $\Proj(R)$ and interpret the
  assumption ``$A$ good determinantal'' by ``$A$ standard determinantal
  satisfying $\depth_{I_{t-1}(\cA)A}A \ge 1$'' (\cite{KMMNP}, Prop.\! 3.2).
  Then the proof of Theorem~\ref{compthm} works (we need Remark~\ref{comprem})
  since we have $ \Hom_R(M,M) \simeq A$ by Remark~\ref{dimWba2}. Using
  Remark~\ref{Amodulerem} we get that Theorem~\ref{Amodulethm3},
  Theorem~\ref{Amodulethm5} and Theorem~\ref{normal}(i) are valid in this
  generality while it for the corollaries and Theorem~\ref{normal} (ii)
  suffices to suppose that $\Proj(R)$ is a smooth ACM-scheme (in the case $c >
  2$, see the next theorem for $c = 2$).
  Note that the assumption $k=\overline{k}$ allows us to keep the definition ${
    W(\underline{b};\underline{a})}$ as a certain locus in $\Hi ^p(\PP^{n})$.
 \end{remark}

 Finally we will illustrate the results mentioned in the last remark to see
 that, in addition to reproving and generalizing Ellingsrud's codimension $2$
 result (\cite{elli}) a little, we can enlighten the differences between the
 cases $c=2$ and $c > 2$. Indeed the main ingredient is that if $c=2$ and
 $X=\Proj(A)$ is standard determinantal in an ACM scheme $Y=\Proj(R)$, then $M
 \simeq K_A(s)$ for some integer $s$ where $K_A$ is the canonical module of
 $A$ (cf. the line before \eqref{ell}). It follows that we do not need the
 results of section 4 at all to conclude that $ _0\! \Ext_A^i(M,M) = 0$ for
 $i>0$ because this is well known. Moreover in section 5 we needed the weak
 assumption $\depth_{I_{t-1}(\cA)A}A \ge 1$ to get $\Hom_A(M,M) \simeq A$
 which was central in \eqref{specseq}\! -\eqref{specseq2} and hence in the
 proof of Theorem~\ref{compthm}. Now this isomorphism always holds, again by
 $M \simeq K_A(s)$, and we get $ {\rm Def}_{M/R} \simeq {\rm Def}_{A/R}$
 without requiring $\depth_{I_{t-1}(\cA)A}A \ge 1$. These functor are formally
 smooth (Theorem~\ref{modulethm}, Remarks~\ref{modulerem} and \ref{dimWba2})
 and we deduce the theorem below where we interpret $ \Hi (Y)$ as $ \Hi
 ^p(Y)$ (resp. $\GradAlg(H)$) if $\dim X \ge 1$ (resp. $\dim X=0$) as in
 Theorem~\ref{Amodulethm3}. Notice that we now deal with standard
 determinantal schemes $X$ of codimension $2$ {\em in $Y =\Proj(R)$} (they are
 usually not determinantal schemes in $ \PP^{n}$). With $\underline{b}$,
 $\underline{a}$ as in \eqref{minpres} and $X \in
 W_s(\underline{b};\underline{a})$ we get
 \begin{theorem} \label{elling} Let $Y=\Proj(R) \subset \PP_k^{n}$ be an ACM
   scheme where $k$ is any field and let $X=\Proj(A) \subset Y$,
   $A=R/I_{t}(\cA)$, be any standard determinantal scheme of codimension $2$ in
   $Y$. Then $\Hi (Y)$ is smooth at $(X)$ and \ $\dim_{(X)} \Hi (Y) =
   \lambda(R)_2$ where
 $$ \lambda(R)_2:= \sum_{i,j} \dim R_{(a_i-b_j)}
 + \sum_{i,j} \dim R_{(b_j-a_i)} - \sum _{i,j} \dim R_{(a_i-a_j)}- \sum _{i,j}
 \dim R_{(b_i-b_j)} + 1 . $$ Moreover every deformation of $X$ comes from
 deforming $\cA$. In particular if $k=\overline{k}$, then $\Hi (Y)$ is smooth
 along $W_s(\underline{b};\underline{a})$ and the closure $\overline{
   W_s(\underline{b};\underline{a})}$ in $\Hi (Y)$ is an irreducible component
 of dimension $ \lambda(R)_2$.
\end{theorem}
Indeed there are no singular points $(X)$ of  $\Hi (Y)$, $(X) \in
W_s(\underline{b};\underline{a})$ while singular points of  $\Hi (Y)$ for
$c>2$ at $(X) \in W_s(\underline{b};\underline{a})$ are quite common
(see \cite{MDPi} and Rem.\! 3.6 of \cite{KM09}).

\end{document}